%% file: main.tex
\documentclass{article}
\usepackage[a4paper,margin=1in,footskip=0.25in]{geometry}
\usepackage{graphicx} 
\usepackage{amsmath}
\usepackage{amssymb}
\usepackage{authblk}
\usepackage{hyperref}
\usepackage{mathtools}
\usepackage{stmaryrd}
\usepackage{xfrac}
\usepackage{natbib}
\usepackage{siunitx}
\usepackage{color}
\usepackage{tcolorbox}
\usepackage{subcaption}
\usepackage{tikz}
\usepackage{multirow}
\usepackage{tikz-3dplot}
 \usetikzlibrary{calc}
 \usetikzlibrary{patterns}

\newcommand{\twodots}{\mathinner {\ldotp \ldotp}}

\def\correspondingauthor{\footnote{Corresponding author: louis.libat2@univ-eiffel.fr}}

\title{A Cartesian Cut-Cell Two-Fluid Method for Two-Phase Diffusion Problems}
\author[1]{Louis Libat \correspondingauthor}
\author[1]{Can Selçuk}
\author[1]{Eric Chénier}
\author[1]{Vincent Le Chenadec}
\affil[1]{MSME, Université Gustave Eiffel, UMR CNRS 8208, Marne-la-Vallée, 77454, France}
\date{\today}

\begin{document}

\maketitle

\begin{abstract}
We present a Cartesian cut-cell finite-volume method for sharp-interface
two-phase diffusion problems in static geometries. The formulation
follows a two-fluid approach: independent diffusion equations are discretized
in each phase on a fixed Cartesian grid, while the phases are coupled
through embedded interface conditions enforcing continuity of diffusive flux and
a general jump law.
Cut cells are treated by integrating the governing equations over
phase-restricted control volumes and surfaces, yielding discrete divergence and
gradient operators that are locally conservative within each phase.
Interface coupling is achieved by introducing a small set of interfacial
unknowns per cut cell on the embedded boundary; the
resulting algebraic system involves only bulk and interfacial averages. A key feature of the method is the use of a reduced set of geometric information based solely on low-order moments (trimmed volumes, apertures and interface measures/centroids),
allowing robust implementation without constructing explicitly cut-cell polytopes.
The method supports steady (Poisson) and unsteady (diffusion) regimes and
incorporates Dirichlet, Neumann, Robin boundary conditions and general jumps. We validate the scheme on one-, two- and three-dimensional single-phase and two-phase
benchmarks, including curved embedded boundaries, Robin conditions and strong
property/jump contrasts. The results demonstrate a superlinear convergence
behavior, sharp enforcement of interfacial laws and excellent conservation
properties. Extensions to moving interfaces and Stefan-type free-boundary
problems are natural perspectives of this framework.
\end{abstract}

\section{Introduction}

Diffusive transport across material interfaces is central to a wide range of
heat- and mass-transfer phenomena, including conjugate heat transfer between
solids and fluids, interphase species transfer, dissolution/precipitation in
porous media and phase-change processes. At the continuum level, these
problems are naturally formulated as two-phase diffusion models in which
a scalar field (e.g.\ temperature or concentration) satisfies diffusion
equations with phase-dependent coefficients and is coupled across a sharp
interface through jump relations and flux balance
\citep{bird_transport_2007,cussler_diffusion_2009,carslaw_conduction_1980,Ishii2011}.
Depending on the application, the interface may enforce continuity of the
scalar (conjugate heat transfer), a weighted jump such as Henry's law
(partitioning in mass transfer) or a finite interfacial resistance such as
Kapitza resistance \citep{Sander2023,pollack_kapitza_1969}. These interfacial
relations must be satisfied accurately to predict transfer rates, which are
often limited by diffusive fluxes localized near the interface.

From a numerical standpoint, sharp-interface diffusion on complex geometries is
challenging whenever the computational mesh does not conform to the interface.
Body-fitted (ALE-type) approaches can represent the interface explicitly but
require mesh generation and mesh motion, with the attendant complexity and
robustness issues in three dimensions \citep{hirt_arbitrary_1974}. On fixed
Cartesian grids, a broad class of non-conforming strategies has therefore been
developed. Immersed-boundary and ghost/immersed-interface-type methods modify
stencils near the interface or introduce forcing terms so that boundary and
interface conditions are satisfied approximately, often with excellent accuracy
for smooth solutions but with less direct control over strict local
conservation \citep{peskin_flow_1972,gibou_second-order-accurate_2002,gabbard_high-order_2024, verzicco_immersed_2023}.
Interface-capturing methods (VOF or level-set approaches) \citep{hirt_volume_1981,osher_fronts_1988,popinet_accurate_2009} coupled with a one-fluid formulation offer flexibility for complex topological changes; however, enforcing
two-phase diffusion jump conditions sharply on a fixed grid typically
requires additional reconstruction and coupling machinery (and, in practice,
often relies on one-fluid effective-property closures in interfacial cells).
Hybrid combinations have recently been explored for solidification/melt
problems \citep{limare_hybrid_2023} and XFEM/level-set formulations have also
been used for coupled Stefan/mass-transport settings \citep{li_numerical_2019},
highlighting continued interest in sharp and accurate interfacial transfer.

A particularly attractive class of methods for diffusion on Cartesian grids is
the family of embedded-boundary (cut-cell) finite-volume schemes, which
retain a conservative flux balance by trimming control volumes intersected by
the interface. Seminal works include the embedded-boundary Poisson solver of
\citep{johansen_cartesian_1998}, the Cartesian-grid advection-diffusion method
of \citep{calhoun_cartesian_2000} and the three-dimensional heat/Poisson
embedded-boundary method of \citep{schwartz_cartesian_2006}. These approaches
demonstrate that cut-cell finite volumes can combine geometric flexibility
with conservation and sharp boundary enforcement. Extensions and related
developments have been pursued for higher-order finite-volume discretizations
and locally refined grids \citep{mccorquodale_high-order_2011}, as well as for
viscous incompressible flows in complex geometries on staggered meshes
\citep{verstappen_symmetry-preserving_2004,rodriguez_conservative_2022}.
Despite this progress, comparatively fewer works provide a two-fluid,
conservative cut-cell formulation tailored to two-phase
diffusion where the diffusion equations are solved independently in each
phase and general jump laws are enforced sharply at the embedded
interface, without resorting to smeared effective properties in mixed cells.

In this paper, we introduce a conservative Cartesian cut-cell finite-volume
framework for two-phase diffusion problems in static geometries. The
interface is fixed in time, while the scalar fields may be obey to a steady
(Poisson) or unsteady diffusion equation. The method follows a two-fluid philosophy:
each phase carries its own bulk unknowns and cut cells contain
two phase-restricted control volumes. To enforce interfacial jump conditions
in a sharp and locally conservative manner, we introduce, in addition to bulk
cell averages, a small set of interfacial unknowns per cut cell representing interface averages. The discrete
operators are derived by consistent application of Gauss' theorem on the
trimmed volumes and faces, leading to conservative divergence and gradient
operators whose interface contributions appear only through geometric weights
and local coupling terms.

A key design goal is practical robustness: rather than explicitly constructing
trimmed polytopes for every intersection pattern, the method is written
entirely in terms of a reduced set of geometric moments (phase volumes,
face apertures, interface measures and centroids). These moments can be
computed using dedicated quadrature methods for implicitly defined interfaces
\citep{saye_high-order_2015}, VOF-type geometric tools \citep{chierici_optimized_2022},
or general-purpose geometry engines in two dimensions \citep{santilli_geos_2023}.
This moment-based viewpoint keeps the algebraic structure of the discretization
independent of the geometric complexity, while preserving strict finite-volume
conservation.

The main contributions of this work are:
\begin{itemize}
  \item a conservative two-fluid cut-cell finite-volume discretization
  for steady and unsteady diffusion with phase-dependent coefficients and
  general sharp jump laws at an embedded interface;
  \item a unified enforcement mechanism for Dirichlet, Neumann and Robin
  conditions as well as weighted flux and scalar jumps that relies only on local geometric
  moments and nearest-neighbor couplings;
  \item a validation on one-, two- and three-dimensional benchmarks,
  mono- and two-phase, demonstrating a super-linear
  accuracy behavior together with
  conservation properties at the discrete level.
\end{itemize}

The paper is organized as follows. Sec.~\ref{sec:continuum} introduces the
two-phase continuum model and the interfacial jump conditions considered.
Sec.~\ref{sec:discrete} presents the embedded-boundary geometric setting and
the reduced moment description used throughout. Sec.~\ref{sec:cutcell}
derives the cut-cell finite-volume operators and the coupled discrete system.
Numerical results are then reported in Sec.~\ref{sec:numstatic} to assess accuracy, robustness and
conservation and Sec.~\ref{sec:conclusion} concludes and outlines extensions.
\section{Continuum modeling}
\label{sec:continuum}

\begin{figure}
\centering
\includegraphics[width=\linewidth]{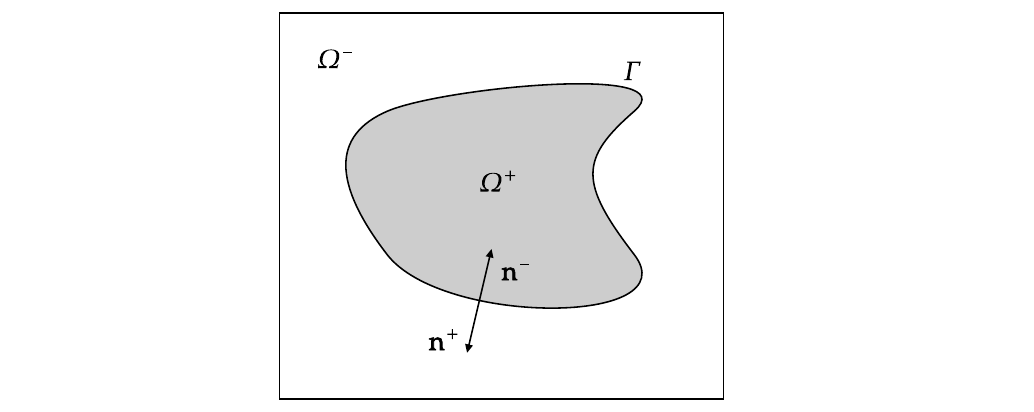}
\caption{Fixed domain \(\Omega\) with interface \(\Gamma\), partitioning into \(\Omega ^ -\) and \(\Omega ^ +\). Normals \(\mathbf n ^ \pm\) point out of their respective phases.}
\label{fig:open-sine-interface}
\end{figure}

Transfers between two immiscible phases are considered in an open domain $\Omega \subset \mathbb{R} ^ d$ ($d \in \{1,2,3\}$). The domain $\Omega$ is partitioned into three subsets, a representative configuration of which is depicted in Fig.~\ref{fig:open-sine-interface}: the first two correspond to the physical regions occupied by the light, $\Omega ^ -$ and dark, $\Omega ^ +$, phases, while the last, $\Gamma$, corresponds to the sharp interface that separates the two phases. The immiscibility condition then reads,
\begin{equation}
\Omega ^ -  \cap \Omega ^ +  = \emptyset.
\label{eq:immiscibility}
\end{equation}
and the saturation condition,
\begin{equation}
\Omega ^ -  \cup \Gamma  \cup \Omega ^ +  = \Omega.
\label{eq:saturation}
\end{equation}
The superscripts ``$-$'' and ``$+$'' denote quantities in the light and dark
phases, respectively and $\mathbf n^\pm$ are the unit normals on $\Gamma$
pointing out of $\Omega^\pm$. We choose the light phase as reference and set
$\mathbf n := \mathbf n^-$, so that interfacial operators are defined with
respect to $\mathbf n$. We employ the symbol $:=$ to mean ``is defined as''.

We introduce a generic scalar field $\phi$ (e.g.\ temperature $T$ or
mass fraction $x$). The following describes the bulk transport model and the
interfacial matching conditions used to close the two-phase diffusion problem.

\subsection{Bulk  transport equations}

In each phase $\Omega ^ \pm$, $\phi ^ \pm$ satisfies a balance equation of the form,
\begin{equation}
C^\pm \,\frac{\partial \phi^\pm}{\partial t} + \nabla \cdot \mathbf q^\pm = r^\pm
\quad t>0,\ \mathbf x \in \Omega^\pm,
\label{eq:bulk-balance}
\end{equation}
where the flux is assumed to be purely diffusive and closed by a constitutive equation such as Fourier's or Fick's law depending on whether the flow of heat or chemical species is considered,
\begin{equation}
\mathbf q^\pm := - K^\pm \nabla \phi^\pm,
\label{eq:constitutive}
\end{equation}
where $K^\pm$ is the diffusive mobility. The capacity $C^\pm$ collects the storage term (e.g., $C^\pm=\rho^\pm c_p^\pm$ and $K^\pm=k^\pm$ for heat; $C^\pm=\rho^\pm$ and $K^\pm=D^\pm$ for species). In the context of this given article, the capacities $C^\pm$ are assumed constant to maintain the linearity of the partial differential equations (PDEs). Variables capacities can also be treated but they introduce non-linearities which go beyond the scope or the current contribution. The symbol $r ^ \pm$ represents any volumetric source or sink term.

The steady diffusion (Poisson) problem is recovered by neglecting transient
storage effects, i.e.\ by setting $\partial_t \phi^\pm = 0$ in
Eq.~\eqref{eq:bulk-balance}. In this case, the governing equations reduce to
\begin{equation}
\nabla \cdot \mathbf q^\pm = r^\pm,
\qquad \mathbf x \in \Omega^\pm,
\end{equation}
supplemented by the same interfacial and boundary conditions as in the unsteady
case.

\subsection{Interfacial conditions} \label{sec:interfacial-conditions}

Let $\lambda$ denote a scalar weight (typically constant) defined on the light domain ($\lambda \colon \Omega ^ - \to \mathbb{R}$). At any interfacial location $\mathbf x \in \Gamma $, we can then define the following weighted jump relative to the light phase,
\begin{equation}
\left \llbracket \phi \right \rrbracket ^ \lambda \left ( t, \mathbf x \right ) := \lim _ {\epsilon \to 0 ^ +}
\left [ \phi ^ + \left ( t, \mathbf x - \epsilon \mathbf n ^ + \right ) - \left ( \lambda \phi ^ - \right ) \left ( t, \mathbf x - \epsilon \mathbf n ^ - \right ) \right ],
\label{eq:weighted-jump}
\end{equation}
where $\mathbf{n} ^ \pm$ are evaluated at $\mathbf{x}$. As will be noted below, the matching conditions used to close the equations will often fallback to the case where $\lambda$ is identically one on $\Omega ^ -$. Consequently, the following shorthand notation is introduced to denote the usual jump operator,
$$
\left \llbracket \phi \right \rrbracket := \left \llbracket \phi \right \rrbracket ^ 1.
$$

The first type of interfacial condition arises from applying conservation principles at the interface. The conservation principle applied to the surface states that for $\phi$ to be conserved (and in the absence of line fluxes, surface accumulation and/or source/sink~\cite{Ishii2011}), the fluxes normal to the interface  cancel out and this reads,
\begin{equation}
\left \llbracket \mathbf q \cdot \mathbf n \right \rrbracket \left ( t, \mathbf x \right ) = 0,
\quad t>0,\ \mathbf x\in\Gamma,
\label{eq:surface-balance}
\end{equation}
In the case of conjugate heat transfer, Eq.~\eqref{eq:surface-balance} corresponds to
continuity of the normal heat flux.

The balance equation Eq.~\eqref{eq:surface-balance} is not sufficient to close the problem Eq.~\eqref{eq:bulk-balance}. The flux balance \eqref{eq:surface-balance} must be complemented by one additional interfacial relation. In the context of conjugate heat transfer, for example, the jump of the dependent variable is typically imposed,
\begin{equation}
\left \llbracket \phi \right \rrbracket \left ( t, \mathbf x \right ) = f, \quad t > 0, \quad \mathbf x \in \Gamma,
\label{eq:continuity}
\end{equation}
whereas in conjugate mass transfer the concentration on either side may be related by Henry's law,
\begin{equation}
\left \llbracket \phi \right \rrbracket ^ {\operatorname{He}} \left ( t, \mathbf x \right ) = f, \quad t > 0, \quad \mathbf x \in \Gamma
\label{eq:henry}
\end{equation}
where $\operatorname{He}$ represents Henry's law constant, which tends to $\operatorname{He}\to 1$ when the solubilities match across the interface~\cite{cussler_diffusion_2009,Sander2023}. Also note that both Eqs.~\eqref{eq:continuity}-\eqref{eq:henry} include discontinuities (non-zero right-hand side $f$) to accommodate source terms \cite{pollack_kapitza_1969}.

The proposed framework can also accommodate more complex boundary conditions, such as nonlinearities (rate laws encountered in heterogeneous kinetics) and mixed boundary conditions (e.g. reservoir models with effective transfer coefficients). These extensions are not discussed here for brevity. As highlighted below, this is possible because the intrinsic diffusion equations on each domain are solved simultaneously, as opposed to a single transport equation for an effective fluid that occupies the union of both domains as in the one-fluid approach. This alleviates the modeling issue of effective properties in interfacial cells and enables the use of different models in each phase (convection, diffusion models...).

\section{Discrete geometric representation} \label{sec:discrete}

The continuum models introduced in Sec.~\ref{sec:continuum} are discretized on a
Cartesian grid. Focusing first on the representation of interfaces, this
section introduces the discrete geometric quantities required by the cut-cell
finite-volume operators.

A wide variety of techniques exist to represent interfaces. For simple objects
(cuboids, cylinders, spheres, etc.), constructive solid geometry (CSG) is a
convenient option~\cite{laidlaw_constructive_1986}. More complex geometries can
be approximated using simplicial meshes or splines (e.g.\ Bézier curves,
NURBS) or described implicitly by a level-set function~\cite{osher_fronts_1988},
in which a signed-distance field is either specified analytically or sampled on the background mesh and
interpolated between nodes. The choice of interface representation is closely tied to the numerical strategy used to solve the surrounding fields. The Arbitrary Lagrangian–Eulerian (ALE) formulation, for instance, is well suited when the interface is discretized with segments (in $d=2$) or triangles (in $d=3$) but it requires moving unstructured meshes. On fixed structured meshes (predominantly Cartesian), the flows on either side of the interface can instead be computed using the cut-cell method, a finite-volume approach typically classified as an embedded boundary method. This technique has been applied to the solution of the incompressible (single-phase) Navier-Stokes equations in complex geometries on staggered meshes, first by Veldman \& Verstappen \cite{verstappen_symmetry-preserving_2004} and later by Botella \cite{cheny_ls-stag_2010}. In both cases, control volumes are explicitly constructed by determining the intersection of the interface with the Cartesian mesh. Once this construction is performed, differentiation operators are derived using the finite-volume methodology.

A major practical difficulty of explicit cut-cell constructions is the large
number of topologically distinct intersection patterns. In a $d$-dimensional
Cartesian cell (a hyperrectangle), the interface can in principle generate
$2^{2^d}$ distinct in/out sign configurations at the $2^d$ vertices. For $d=3$,
this already yields $256$ cases. Adding time as a fourth dimension yields a 4D hyperrectangle with $2^{4}=16$ vertices and $2^{2^{4}}=65536$ possible sign configurations. Although symmetries can be exploited to reduce this number, the resulting independent cases are still considerable.

For this reason, the proposed method does not rely on an explicit construction of trimmed polyhedra. Instead, it relies on Green-Gauss theorem and quadrature rules to derive the discrete operators. These quadrature rules require only a reduced set of geometric moments, described in Sec.~\ref{sec:moments}. The computation of these moments thus constitutes the main bottleneck of the proposed method. It can be carried out with a handful of computational-geometry techniques and libraries. As shown below, this algorithmic choice is facilitated by the fact that all required moments reduce to intersections of hyperrectangles with the light (or dark) domain. For this purpose, efficient algorithms
exist. For implicitly defined geometries, Saye's high-order
quadrature methods provide accurate surface/volume integration directly from a
level-set description~\cite{saye_high-order_2015}. The VOFI library also offers
high-precision volume-of-fluid integration and related geometric quantities
for implicitly defined interfaces~\cite{chierici_optimized_2022}. In two
dimensions, general-purpose geometry engines such as GEOS/LibGEOS provide robust
polygon operations (intersection, area, centroid) ~\cite{santilli_geos_2023}, which can be useful for handling complex planar geometries. Altogether, such tools substantially simplify the implementation of the proposed method.

\subsection{Geometrical foundations}

We now introduce all of geometric quantities required to assemble the discrete operators in Sec.~\ref{sec:cutcell}. For the sake of brevity, this section focuses on the two-dimensional case ($d = 2$). Three-dimensional extension can readily be performed.

Let $X := \left ( x _ {\sfrac{1}{2}}, x _ {\sfrac{3}{2}}, \ldots, x _ {N+\sfrac{1}{2}} \right ) \in \mathbb{R} ^ {N + 1}$ ($N \in \mathbb{N} ^ *$) and $Y:=\left ( y _ {\sfrac{1}{2}}, y _ {\sfrac{3}{2}}, \ldots, y _ {M+\sfrac{1}{2}} \right ) \in \mathbb{R} ^ {M + 1}$ ($M \in \mathbb{N} ^ *$) denote two vectors of strictly increasing abscissas,
$$
x _ {\sfrac{1}{2}} < x _ {\sfrac{3}{2}} < \ldots < x _ {N + \sfrac{1}{2}} \quad \mathrm{and} \quad y _ {\sfrac{1}{2}} < y _ {\sfrac{3}{2}} < \ldots < y _ {M + \sfrac{1}{2}}.
$$

We also define the following notation,
$$
\Delta ^ 1 _ i := \left ] x _ {i - \sfrac{1}{2}}, x _ {i + \sfrac{1}{2}} \right [ \quad \mathrm{and} \quad \Delta ^ 2 _ j := \left ] y _ {j - \sfrac{1}{2}}, y _ {j + \sfrac{1}{2}} \right [,
$$
where the superscripts indicate the directions and where
$$
\left ] a, b \right [ := \left \{ \xi \in \mathbb{R} \enskip \mathrm{s.t.} \enskip a < \xi < b \right \}, \quad \left ( a, b \right ) \in \mathbb{R} ^ 2
$$
denotes a finite open interval.

The rectilinear mesh is defined using the Cartesian product: in two spatial dimensions,
$$
\Omega _ {i, j} := \Delta ^ 1 _ i \times \Delta ^ 2 _ j, \quad \left ( i, j \right ) \in \left [ 1 \twodots N \right ] \times \left [ 1 \twodots M \right ]
$$
where,
$$
\left [ m \twodots n \right ] := \left \{ i \in \mathbb{Z} \enskip \mathrm{s.t.} \enskip m \le i \le n \right \}, \quad \left ( m, n \right ) \in \mathbb{Z} ^ 2
$$
denotes a set of consecutive integers.

\begin{figure}[h]
  \centering
  \begin{subfigure}[b]{0.48\linewidth}
    \centering
    \includegraphics[width=\linewidth,trim={5cm 0 5cm 0},clip]{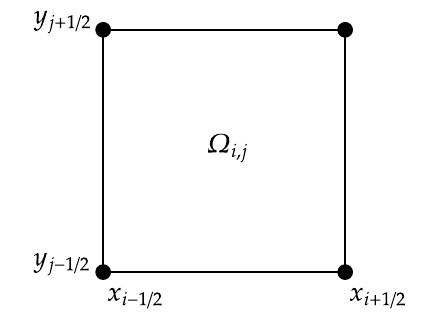}
    \caption{Mesh cell \(\Omega_{i, j}\) with discrete node convention}
    \label{fig:fig2-schema}
  \end{subfigure}
  \quad
  \begin{subfigure}[b]{0.48\linewidth}
    \centering
    \includegraphics[width=\linewidth,trim={5cm 0 5cm 0},clip]{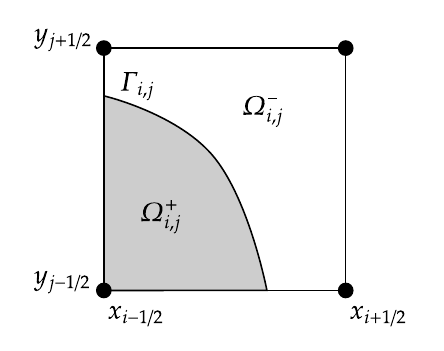}
    \caption{Mixed cell \(\Omega_{i,j}\) with \(\Omega_{i,j}^\pm\) domains and \(\Gamma_{i,j}\)}
    \label{fig:fig3-schema}
  \end{subfigure}
  \caption{Mesh cell notations}
  \label{fig:fig2-3}
\end{figure}

To describe the intrinsic properties of each phase, the intersection of the mesh cells with each phase's domain is introduced as follows: let,
\begin{equation}
\Omega ^ \pm _ {i, j}  := \Omega _ {i, j} \cap \Omega ^ \pm 
\label{eq:cutcell}
\end{equation}
denote the subset of $\Omega _ {i, j}$ occupied by the light/dark phase and
$$
\Gamma _ {i, j} := \Omega _ {i, j} \cap \Gamma
$$
the subset of interface within $\Omega _ {i, j}$. From the immiscibility condition Eq.~\eqref{eq:immiscibility}, it follows that,
\begin{equation} 
\Omega _ {i, j} = \Omega _ {i, j} ^ -  \cup \Gamma _ {i, j}  \cup \Omega _ {i, j} ^ + ,
\label{eq:omega-subset}
\end{equation}
where all the subsets on the right-hand side are disjoint (See Fig \ref{fig:fig3-schema}).

These notations enable the following definitions. A cell $\left ( i, j \right )$ is referred to as purely dark (resp., light) if $\Omega ^ + _ {i, j}  = \Omega _ {i, j}$ (resp., $\Omega ^ - _ {i, j}  = \Omega _ {i, j}$). From the immiscibility condition, it follows that a purely dark cell is void of light medium (and vice versa). A cell $\left ( i, j \right )$ is said to be mixed if $\Gamma _ {i, j}  \ne \emptyset$.

We also introduce the notation,
$$
\Omega ^ {1 \pm} _ j \left (\Delta \right ) := \left ( \Delta \times \Delta ^ 2 _ j \right ) \cap \Omega ^ \pm  \quad \mathrm{and} \quad \Omega ^ {2 \pm} _ i \left (\Delta \right ) := \left ( \Delta ^ 1 _ i \times \Delta \right ) \cap \Omega ^ \pm 
$$
where $\Delta$ denotes a open interval. Note that, $\Omega ^ {1 \pm} _ j \left (\Delta \right )$ and $\Omega ^ {2 \pm} _ i \left (\Delta \right )$ are both subsets of $\mathbb R ^ 2$.

As was previously done for the cells, the intersection of mesh faces with the phases' domains is also introduced as follows. Let
$$
\Sigma ^ {1 \pm} _ j \left (x \right ) := \left \{ y \in \Delta ^ 2 _ j \ \ \mathrm{s.t.} \ \left ( x, y \right ) \in \Omega ^ \pm \right \}, \quad j \in \left [ 1 \twodots M + 1 \right ]
$$
and
$$
\Sigma ^ {2 \pm} _ i \left ( y \right ) := \left \{ x \in \Delta ^ 1 _ i \ \ \mathrm{s.t.} \ \left ( x, y \right ) \in \Omega ^ \pm  \right \}, \quad i \in \left [ 1 \twodots N + 1 \right ].
$$
Note that for all $\left ( x, y \right ) \in \mathbb R ^ 2$, $\Sigma ^ {1 \pm} _ j \left (x \right )$
and $\Sigma ^ {2 \pm} _ i \left ( y \right )$
are both subsets of $\mathbb{R}$.

The description is completed by defining the corresponding faces in $\mathbb{R}^2 $, that is defining both,
$$
\Sigma ^ {1 \pm} _ {i - \sfrac{1}{2}, j}:= \left \{ \left ( x, y \right ) \in \left\{ x _ {i - \sfrac{1}{2}} \right\} \times \Delta ^ 2 _ j \ \ \mathrm{s.t.} \ \left ( x, y \right ) \in \Omega ^ \pm  \right \},
$$
and,
$$
\Sigma ^ {2 \pm} _ {i, j - \sfrac{1}{2}}  := \left \{ \left ( x, y \right ) \in \Delta ^ 1 _ i \times \left\{ y _ {j - \sfrac{1}{2}} \right\} \ \ \mathrm{s.t.} \ \left ( x, y \right ) \in \Omega ^ \pm  \right \},
$$

\begin{figure}[htbp]
  \centering
\includegraphics[width=\linewidth]{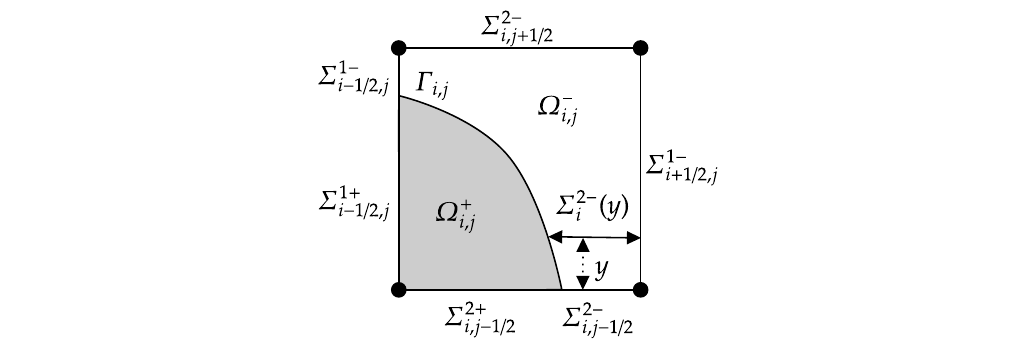}
  \caption{Mesh cell \(\Omega_{i, j}\) with interface \(\Gamma_{i, j}\) separating dark \(\Omega_{i, j}^+\) and light regions \(\Omega_{i, j}^-\). Intersections of mesh faces with phases' domains is also represented}
  \label{fig:mesh-cell}
\end{figure}

The qualifiers employed earlier for the cells (pure and mixed) also apply to these faces.
Finally, the subsets thus defined realize a partition of the boundary (vertices excluded) of  each phase (See Fig \ref{fig:mesh-cell}),
\begin{equation}
\partial \Omega _ {i, j} ^ \pm \simeq \Gamma _ {i, j} \cup \Sigma ^ {1 \pm} _ {i - \sfrac{1}{2}, j} \cup \Sigma ^ {1 \pm} _ {i + \sfrac{1}{2}, j} \cup \Sigma ^ {2 \pm} _ {i, j - \sfrac{1}{2}} \cup \Sigma ^ {2 \pm} _ {i, j + \sfrac{1}{2}}.
\label{eq:partition}
\end{equation}

\subsection{Reduced geometric description} \label{sec:moments}

Let $\Xi$ denote a subset of $\mathbb{R} ^ d$ ($d \in \mathbb{N} ^ *$) and $f$ a function defined over $\Xi$ ($f$ is typically a monomial, as shown below). We denote the integral of $f$ over $\Xi$ with the following shorthand notation,
\begin{equation}
\left \langle \Xi, f \right \rangle := \int _ {\left ( x _ 1 \twodots x _ d \right ) \in \Xi} f \left ( x _ 1, \ldots, x _ d \right ) \, \mathrm{d} x _ 1 \ldots \mathrm{d} x _ d
\label{eq:integral}
\end{equation}
which will be employed over cells, faces and interface segments. (Here $x _ 1$, $\ldots$, $x _ d$ denote Cartesian coordinates).

In the computer‐graphics literature these same quantities are often called geometric moments \cite{mirtich_fast_1996}. These same moments are employed in various cut-cell discretizations to compute irregular cell volumes and face apertures needed for flux evaluation~\cite{remmerswaal_sharp_2022}. We also note that the notation introduced in Eq.~\eqref{eq:integral} is sufficient, since all required geometric information is expressed as integrals over the intersection of hyperrectangles with phase domains. 


\subsubsection{Primary geometric moments}

As mentioned above, only low-order moments are required to assemble the discrete operators. This section therefore introduces shorthand notations for such moments. Let us first denote the volumes,
$$
V _ {i, j} := \left \langle \Omega _ {i, j}, 1 \right  \rangle  =\left ( x _ {i + \sfrac{1}{2}} - x _ {i - \sfrac{1}{2}} \right ) \left ( y _ {j + \sfrac{1}{2}} - y _ {j - \sfrac{1}{2}} \right )
$$
and
\begin{equation}
    V ^ \pm _ {i, j} := \left \langle \Omega ^ \pm _ {i, j} , 1 \right \rangle.
\label{eq:volume_cut}
\end{equation}

\begin{figure}[!h]
  \centering
\includegraphics[width=\linewidth]{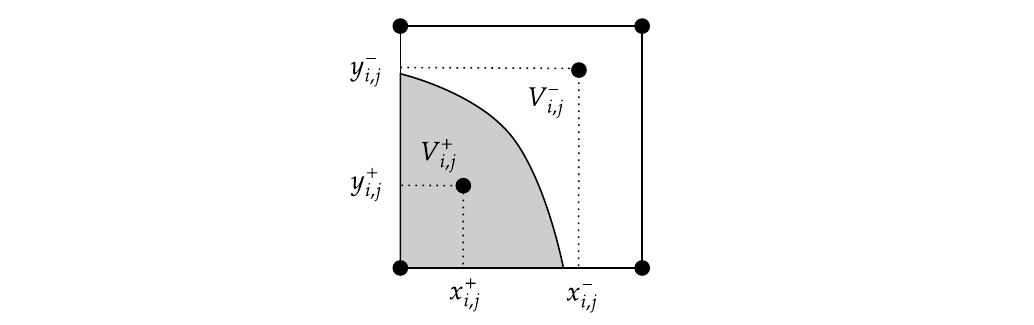}
  \caption{Cut cell with the domains \(\Omega_{i, j}^+\) (dark) and \(\Omega_{i, j}^-\) (light), with bulk centroids \((x_{i, j}^\pm,y_{i, j}^\pm)\) and their associated \(V_{i, j}^\pm\)}
  \label{fig:fig5}
\end{figure}

The centroids along $x$ are defined as (See Fig \ref{fig:fig5}),
$$
x _ i := \frac{\left \langle \Omega _ {i, j}, x \right \rangle}{\left \langle \Omega _ {i, j}, 1 \right \rangle}  = \frac{x _ {i - \sfrac{1}{2}} + x _ {i + \sfrac{1}{2}}}{2} 
$$
and
\begin{equation}
x ^ \pm _ {i, j}  := \begin{cases}
\left \langle \Omega ^ \pm _ {i, j}, x \right \rangle / \left \langle \Omega ^ \pm _ {i, j} , 1 \right \rangle, & \mathrm{if} \enskip \Omega ^ \pm _ {i, j}  \ne \emptyset, \\
x _ i, & \mathrm{otherwise}.
\end{cases}
\label{eq:bulk-center-x}
\end{equation}
Likewise,
$$
y _ j := \frac{\left \langle \Omega _ {i, j}, y \right \rangle}{\left \langle \Omega _ {i, j}, 1 \right \rangle} \left ( = \frac{y _ {j - \sfrac{1}{2}} + y _ {j + \sfrac{1}{2}}}{2} \right )
$$
and
\begin{equation}
y ^ \pm _ {i, j}  := \begin{cases}
\left \langle \Omega ^ \pm _ {i, j} , y \right \rangle / \left \langle \Omega ^ \pm _ {i, j} , 1 \right \rangle, & \mathrm{if} \enskip \Omega ^ \pm _ {i, j}   \ne \emptyset, \\
y _ j, & \mathrm{otherwise}.
\end{cases}
\label{eq:bulk-center-y}
\end{equation}
The definitions Eq.~\eqref{eq:bulk-center-x} and Eq.~\eqref{eq:bulk-center-y} remain valid for empty cells, for which either $V _ {i, j} ^ +$ or $V _ {i, j} ^ -$ vanishes. In such instances, this enables a consistent definitions of the secondary geometric moments, defined in the next section.

This construction also guarantees that the centroids lie within the convex hull of the $\Omega_{i,j}$, i.e.
\begin{equation}
x_{i-\sfrac12}\le x^\pm_{i,j}\le x_{i+\sfrac12},
\qquad
y_{j-\sfrac12}\le y^\pm_{i,j}\le y_{j+\sfrac12}.
\label{eq:convex}
\end{equation}

\begin{figure}[h]
\centering
\includegraphics[width=\linewidth]{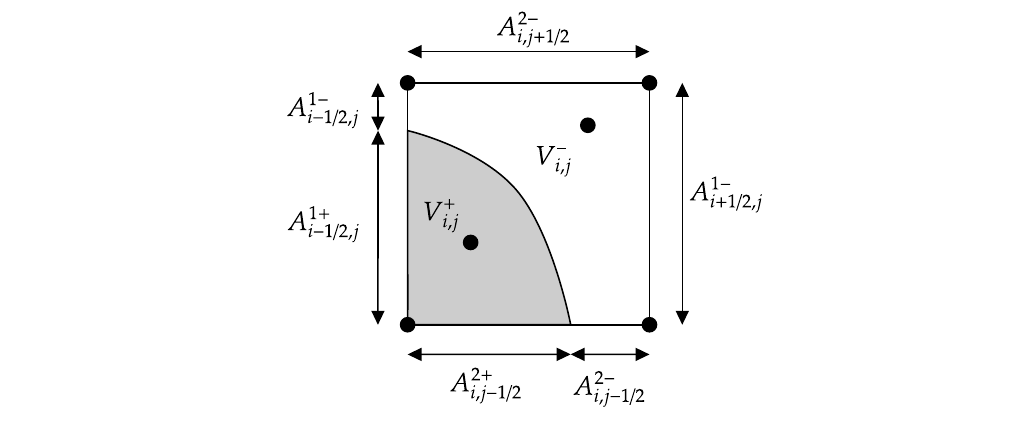}
\caption{Cut cell \(\Omega_{i, j}\) with the domains \(\Omega_{i, j}^+\) (dark) and \(\Omega_{i, j}^-\) (light), with face‐areas \(A^{1\pm}_{i\pm\frac12,j}\), \(A^{2\pm}_{i,j\pm\frac12}\).}
\label{fig:fig6}
\end{figure}

We also introduce the following face-centered quantity (See Fig \ref{fig:fig6}),
\begin{equation}   
A^{1\pm}_{i-\sfrac12,j} := \langle \Sigma^{1\pm}_{i-\sfrac12,j},1\rangle,
\qquad
A^{2\pm}_{i,j-\sfrac12} := \langle \Sigma^{2\pm}_{i,j-\sfrac12},1\rangle.
\label{eq:AxAy}
\end{equation}

Using both the rectilinearity of the grid and the saturation and immiscibility conditions Eqs.~\eqref{eq:immiscibility} and~\eqref{eq:saturation},
$$
V _ {i, j} = V _ {i, j} ^ - + V _ {i, j} ^ +,
$$
$$
x _ i V _ {i, j} = x ^ - _ {i, j} V _ {i, j} ^ - + x ^ + _ {i, j} V _ {i, j} ^ + \quad \mathrm{and} \quad y _ j V _ {i, j} = y ^ - _ {i, j} V _ {i, j} ^ - + y ^ + _ {i, j} V _ {i, j} ^ +
$$
and finally,
$$
A ^ 1 _ j = A ^ {1 -} _ {i - \sfrac{1}{2}, j} + A ^ {1 +} _ {i - \sfrac{1}{2}, j}
\quad \mathrm{and} \quad
A ^ 2 _ i = A ^ {2 -} _ {i, j - \sfrac{1}{2}} + A ^ {2 +} _ {i, j - \sfrac{1}{2}}.
$$
where $A ^ 1 _ j := y _ {j + \sfrac{1}{2}} - y _ {j - \sfrac{1}{2}}$ and $A ^ 2 _ i := x _ {i + \sfrac{1}{2}} - x _ {i - \sfrac{1}{2}}$.

\subsubsection{Secondary geometric moments}

The geometric moments defined below are referred to as ``secondary'' simply because their definition depends on the previously introduced ``primary geometric moment''.

\begin{figure}[h]
\centering
\includegraphics[width=\linewidth]{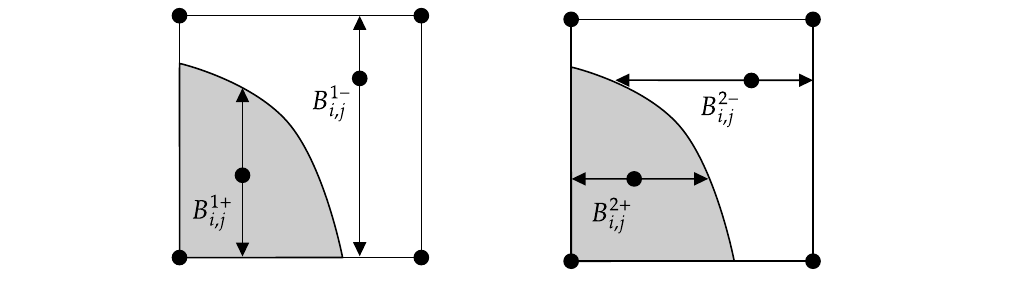}
\caption{Cut cell \(\Omega_{i, j}\) with the domains \(\Omega_{i, j}^+\) (dark) and \(\Omega_{i, j}^-\) (light), with centroid‐areas \(B^{1\pm}_{i,j}\), \(B^{2\pm}_{i,j}\)}
\label{fig:fig7}
\end{figure}

We define the following cell-centered averages (See Fig \ref{fig:fig7}),
\begin{equation}  
B^{1\pm}_{i,j} := \left\langle \Sigma^{1\pm}_j(x^\pm_{i,j}),1\right\rangle,
\qquad
B^{2\pm}_{i,j} := \left\langle \Sigma^{2\pm}_i(y^\pm_{i,j}),1\right\rangle .
\label{eq:BxBy}
\end{equation}

Let us also define the $x$-face-centered volume as,
\begin{equation}
W^{1\pm}_{i-\sfrac12,j}
:= \left\langle
\Omega^{1\pm}_j\!\left(\left]x^\pm_{i-1,j},\,x^\pm_{i,j}\right[\right),
1\right\rangle .
\label{eq:second-volume}
\end{equation}
Using the linearity of the integral, the latter moment can also be computed as,
\[
W^{1\pm}_{i-\sfrac12,j}
=
\left\langle
\Omega^{1\pm}_j\!\left(\left]x^\pm_{i-1,j},\,x_{i-\sfrac12}\right[\right),
1\right\rangle
+
\left\langle
\Omega^{1\pm}_j\!\left(\left]x_{i-\sfrac12},\,x^\pm_{i,j}\right[\right),
1\right\rangle,
\]
Using Eq.~\eqref{eq:convex}, we note that $x ^ \pm _ {i - 1, j} \left ( t \right ) < x _ {i - \sfrac{1}{2}} < x ^ \pm _ {i, j} \left ( t \right )$ and hence the right-hand side is the sum of two positive quantities. Likewise, along $y$-aligned faces we define,
\[
W^{2\pm}_{i,j-\sfrac12}
:= \left\langle
\Omega^{2\pm}_i\!\left(\left]y^\pm_{i,j-1},\,y^\pm_{i,j}\right[\right),
1\right\rangle,
\]
which can also be rewritten as the sum of two positive quantities as done in Eq.~\eqref{eq:second-volume}.

\begin{figure}[h]
\centering
\includegraphics[width=\linewidth]{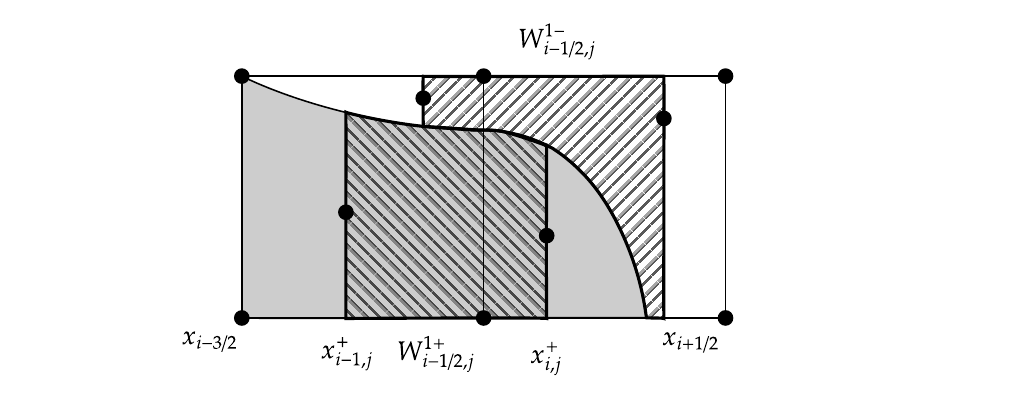}
\caption{Two adjacent cut cells.  The right striped region \(W^{1-}_{i-\frac12,j}\) and the left striped region \(W^{1+}_{i-\frac12,j}\) are shown in the cells.}
\label{fig:fig8}
\end{figure}

A summary table, located in the appendix \ref{tab:geom_and_sd}, summarizes, for two spatial dimensions, the geometric moments used in the discrete formulation.

\section{Two-phase cut-cell method for static domains} \label{sec:cutcell}

We now turn our attention to the discretization of the two scalar transport equations Eq.~\eqref{eq:bulk-balance} on fixed subdomains i.e.\ $\Omega^\pm$. The formulation is based on the single-phase cut-cell method developed in~\cite{rodriguez_conservative_2022}. The novelty
here lies in its extension to two-phase diffusion with interfacial matching
conditions. This extension proceeds as follows: bulk and interfacial variables are introduced in each domain, with the exact number of unknowns depending on the interface topology and the Cartesian mesh. The key question then becomes how to impose discrete constraints that uniquely
determine these unknowns.

Let us introduce a sequence of strictly increasing time instants,
$$
t _ 0 < t _ 1 < \ldots < t _ n < \ldots
$$
where $t _ 0$ is set to $0$, consistently with the continuum model (See Sec.~\ref{sec:continuum})

For each mixed cells, the semi-discrete interfacial variables are defined as,
$$
\Phi _ {i, j} ^ {\gamma \pm} \left ( t \right ) := \frac{\int _ {\Gamma _ {i, j}} \phi^\pm \left ( t , \mathbf x\right ) \, \mathrm{d} S}{\int _ {\Gamma _ {i, j}} \, \mathrm{d} S},
$$
and their discrete-in-time counterparts
\begin{equation}
\Phi _ {n, i, j} ^ {\gamma \pm} := \Phi _ {i, j} ^ {\gamma \pm} \left ( t _ n \right ).
\label{eq:surface-unknown}
\end{equation}
These interfacial unknowns are associated with the two interface conditions (Eqs.~\eqref{eq:surface-balance}-\eqref{eq:continuity} or Eqs.~\eqref{eq:surface-balance}-\eqref{eq:henry} depending on the the problem at hand).

In each active cell, we define the semi-discrete bulk variables as,
$$
\Phi _ {i, j} ^ {\omega \pm} \left ( t \right ) := 
     \left \langle \Omega ^ \pm _ {i, j}, \phi^\pm \left ( t \right ) \right \rangle / V ^ \pm _ {i, j} \quad \mathrm{where} \quad V ^ \pm _ {i, j} \ne 0,
$$
and the discrete bulk variables as,
\begin{equation}
\Phi _ {n, i, j} ^ {\omega \pm} := \Phi _ {i, j} ^ {\omega \pm} \left ( t _ n \right ).
\label{eq:bulk-unknown}
\end{equation}
They are associated with the semi-discrete and discrete bulk equations, respectively.

\subsection{Semi-discrete bulk equations}
\label{sec:discrete_bulk}

As for the primary variables, we represent the Cartesian components of the flux
$\mathbf q^\pm$ by face averages. Let $q^{1\pm} := \mathbf q^\pm\cdot \mathbf e^1$
and $q^{2\pm} := \mathbf q^\pm\cdot \mathbf e^2$, where $\mathbf e^\alpha$ denotes the
$\alpha$-th Cartesian basis vector. On vertical faces, we define
\begin{equation}
Q^{1\pm}_{i-\sfrac12,j}(t)
:= \frac{1}{A^{1\pm}_{i-\sfrac12,j}}
\int_{\Sigma^{1\pm}_{i-\sfrac12,j}} q^{1\pm}(t,\mathbf x)\,\mathrm dS,
\qquad \text{for } A^{1\pm}_{i-\sfrac12,j}\neq 0,
\label{eq:flux-x}
\end{equation}
and similarly on horizontal faces,
\begin{equation}
Q^{2\pm}_{i,j-\sfrac12}(t)
:= \frac{1}{A^{2\pm}_{i,j-\sfrac12}}
\int_{\Sigma^{2\pm}_{i,j-\sfrac12}} q^{2\pm}(t,\mathbf x)\,\mathrm dS,
\qquad \text{for } A^{2\pm}_{i,j-\sfrac12}\neq 0.
\label{eq:flux-y}
\end{equation}

Integrating Eq.~\eqref{eq:bulk-balance} over $\Omega^\pm_{i,j}$ and using the
(divergence) Gauss theorem yields the semi-discrete balance
\begin{equation}
V^\pm_{i,j}\,C^\pm \frac{\mathrm d \Phi^{\omega\pm}_{i,j}}{\mathrm dt}
+ \int_{\partial\Omega^\pm_{i,j}} \mathbf q^\pm(t,\mathbf x)\cdot \mathbf n^\pm \,\mathrm dS
=
\int_{\Omega^\pm_{i,j}} r^\pm(t,\mathbf x)\,\mathrm dV.
\label{eq:bulk-balance-semi}
\end{equation}

\subsubsection{Volume-integrated divergence operator}

Using the partition Eq.~\eqref{eq:partition}, we split the boundary integral into
its mesh-face and interface contributions:
\begin{equation}
\int_{\partial\Omega^\pm_{i,j}} \mathbf q^\pm\cdot \mathbf n^\pm \,\mathrm dS
=
\int_{\partial\Omega^\pm_{i,j}\setminus \Gamma_{i,j}} \mathbf q^\pm\cdot \mathbf n^\pm \,\mathrm dS
+
\int_{\Gamma_{i,j}} \mathbf q^\pm\cdot \mathbf n^\pm \,\mathrm dS.
\label{eq:gauss}
\end{equation}
On mesh-aligned faces, the outward normals are $\pm \mathbf e^1$ and $\pm \mathbf e^2$,
so that the first term can be written in conservative flux-difference form:
\begin{multline}
\int_{\partial\Omega^\pm_{i,j}\setminus \Gamma_{i,j}} \mathbf q^\pm\cdot \mathbf n^\pm \,\mathrm dS
=
A^{1\pm}_{i+\sfrac12,j}\,Q^{1\pm}_{i+\sfrac12,j}(t)
- A^{1\pm}_{i-\sfrac12,j}\,Q^{1\pm}_{i-\sfrac12,j}(t) \\
+ A^{2\pm}_{i,j+\sfrac12}\,Q^{2\pm}_{i,j+\sfrac12}(t)
- A^{2\pm}_{i,j-\sfrac12}\,Q^{2\pm}_{i,j-\sfrac12}(t).
\label{eq:div-regular}
\end{multline}

We introduce directional (volume-integrated) divergence operators acting on generic
apertures $A$ and face fluxes $Q$:
\begin{equation}
\left\{
\begin{aligned}
V^\pm_{i,j}\,\operatorname{div}^{1\omega}_{i,j}(A,Q)
&:= A_{i+\sfrac12,j}\,Q_{i+\sfrac12,j} - A_{i-\sfrac12,j}\,Q_{i-\sfrac12,j},\\
V^\pm_{i,j}\,\operatorname{div}^{2\omega}_{i,j}(A,Q)
&:= A_{i,j+\sfrac12}\,Q_{i,j+\sfrac12} - A_{i,j-\sfrac12}\,Q_{i,j-\sfrac12}.
\end{aligned}
\right.
\label{eq:div-discrete-omega}
\end{equation}

\begin{figure}[!h]
    \centering
    \includegraphics[width=\linewidth]{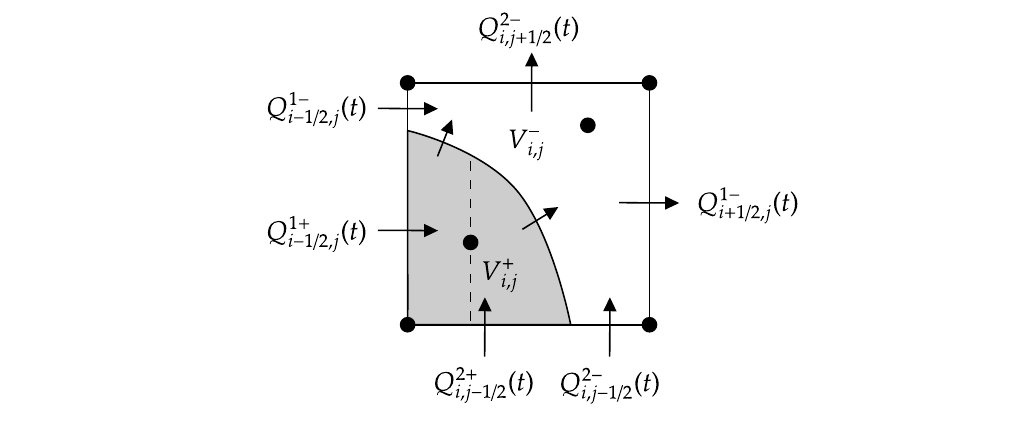}
    \caption{Schematic of the cut-cell divergence balance.}
    \label{fig:cutcell_divergence}
\end{figure}

This step is the foundation of the finite-volume method: each contribution is an
exchange through a well-identified boundary portion. By construction, internal
exchanges cancel under telescopic summation, leaving only net boundary fluxes,
thus ensuring global conservation within each phase.

The last term in Eq.~\eqref{eq:gauss}, representing interfacial exchange, must be
approximated using the available reduced geometric information. The discrete
formula is designed to satisfy the following requirements: (i) minimal stencil,
(ii) exact cancellation when the face fluxes are constant, (iii) at least
first-order accuracy for mesh-aligned planar interfaces and (iv) recovery of
the classical divergence away from interfaces.

We adopt the following approximation:
\begin{multline}
\int_{\Gamma_{i,j}} \mathbf q^\pm(t,\mathbf x)\cdot \mathbf n^\pm \,\mathrm dS
\simeq
\bigl(B^{1\pm}_{i,j}-A^{1\pm}_{i+\sfrac12,j}\bigr)\,Q^{1\pm}_{i+\sfrac12,j}(t)
+\bigl(A^{1\pm}_{i-\sfrac12,j}-B^{1\pm}_{i,j}\bigr)\,Q^{1\pm}_{i-\sfrac12,j}(t) \\
+\bigl(B^{2\pm}_{i,j}-A^{2\pm}_{i,j+\sfrac12}\bigr)\,Q^{2\pm}_{i,j+\sfrac12}(t)
+\bigl(A^{2\pm}_{i,j-\sfrac12}-B^{2\pm}_{i,j}\bigr)\,Q^{2\pm}_{i,j-\sfrac12}(t).
\label{eq:interface-transfer}
\end{multline}
We formalize this by introducing
\begin{equation}
\left\{
\begin{aligned}
V^\pm_{i,j}\,\operatorname{div}^{1\gamma}_{i,j}(A,B,Q)
&:= \bigl(B_{i,j}-A_{i+\sfrac12,j}\bigr)\,Q_{i+\sfrac12,j}
   +\bigl(A_{i-\sfrac12,j}-B_{i,j}\bigr)\,Q_{i-\sfrac12,j},\\
V^\pm_{i,j}\,\operatorname{div}^{2\gamma}_{i,j}(A,B,Q)
&:= \bigl(B_{i,j}-A_{i,j+\sfrac12}\bigr)\,Q_{i,j+\sfrac12}
   +\bigl(A_{i,j-\sfrac12}-B_{i,j}\bigr)\,Q_{i,j-\sfrac12}.
\end{aligned}
\right.
\label{eq:div-discrete-gamma}
\end{equation}

Summing the bulk ($\omega$) and interfacial ($\gamma$) contributions
Eq.~\eqref{eq:div-regular} and Eq.~\eqref{eq:interface-transfer}, the terms proportional
to $A^{\alpha\pm}$ cancel and the divergence can be expressed using $B^{\alpha\pm}$ only:
\[
\left\{
\begin{aligned}
V^\pm_{i,j}\,\operatorname{div}^{1}_{i,j}(B,Q)
&:= B_{i,j}\,\bigl(Q_{i+\sfrac12,j}-Q_{i-\sfrac12,j}\bigr),\\
V^\pm_{i,j}\,\operatorname{div}^{2}_{i,j}(B,Q)
&:= B_{i,j}\,\bigl(Q_{i,j+\sfrac12}-Q_{i,j-\sfrac12}\bigr).
\end{aligned}
\right.
\]
As previously noted, an explicit construction of trimmed control volumes is not
required; the formulation only relies on reduced geometric moments.

\subsubsection{Gradient operator}

\begin{figure}[!h]
    \centering
    \includegraphics[width=\linewidth]{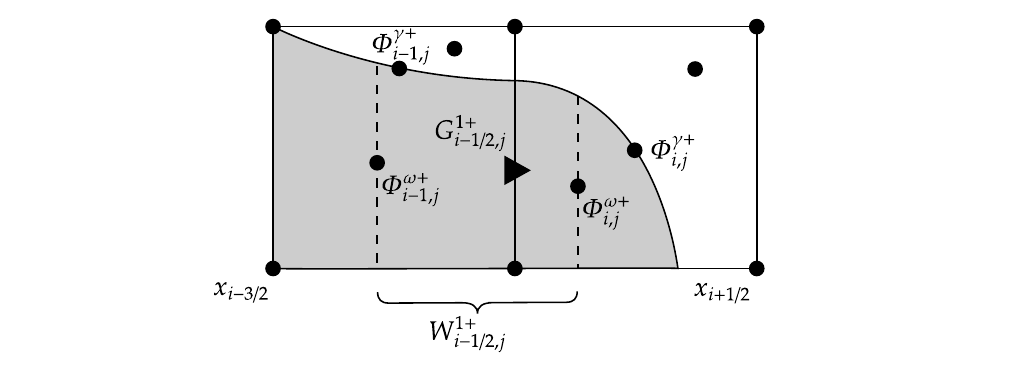}
    \caption{Schematic of the cut-cell gradient operator on the staggered volume
    $W^{1+}_{i-\frac12,j}$.}
    \label{fig:cutcell_gradient}
\end{figure}

The first component of the gradient is represented by a staggered average
\begin{equation}
G^{1\pm}_{i-\sfrac12,j}(t)
:= \frac{1}{W^{1\pm}_{i-\sfrac12,j}}
\left\langle \Omega^{1\pm}_j\!\left(\left]x^\pm_{i-1,j},\,x^\pm_{i,j}\right[\right),
\ \partial_1 \phi^\pm(t,\cdot)\right\rangle,
\qquad \text{for } W^{1\pm}_{i-\sfrac12,j}\neq 0,
\label{eq:grad-x}
\end{equation}
where $\partial_1$ denotes the partial derivative with respect to $x$.
Using the Gauss theorem on the staggered region and the linearity of the
integral, we obtain
\begin{multline}
W^{1\pm}_{i-\sfrac12,j}\,G^{1\pm}_{i-\sfrac12,j}(t)
=
\left\langle \Sigma^{1\pm}_j(x^\pm_{i,j}),\ \phi^\pm(t,\cdot)\right\rangle
-
\left\langle \Sigma^{1\pm}_j(x^\pm_{i-1,j}),\ \phi^\pm(t,\cdot)\right\rangle \\
+ \int_{\left(\left]x^\pm_{i-1,j},\,x^\pm_{i,j}\right[\times \Delta^2_j\right)\cap \Gamma}
\phi^\pm(t,\mathbf x)\,n^{1\pm}(\mathbf x)\,\mathrm dS,
\label{eq:grad_decomposition}
\end{multline}
where $n^{1\pm}$ is the first Cartesian component of the unit normal to $\Gamma$
pointing out of $\Omega^\pm$.

The two mesh-face terms are approximated by
\begin{equation}
\left\langle \Sigma^{1\pm}_j(x^\pm_{i,j}),\ \phi^\pm(t,\cdot)\right\rangle
\simeq B^{1\pm}_{i,j}\,\Phi^{\omega\pm}_{i,j}(t),
\qquad
\left\langle \Sigma^{1\pm}_j(x^\pm_{i-1,j}),\ \phi^\pm(t,\cdot)\right\rangle
\simeq B^{1\pm}_{i-1,j}\,\Phi^{\omega\pm}_{i-1,j}(t),
\label{eq:substitutions}
\end{equation}
which yields the bulk contribution
\[
\operatorname{grad}^{1\omega}_{i-\sfrac12,j}(B,W,\Phi)
:= \frac{B_{i,j}\Phi_{i,j}-B_{i-1,j}\Phi_{i-1,j}}{W_{i-\sfrac12,j}}.
\]

The interfacial contribution (last term in Eq.~\eqref{eq:grad_decomposition}) is
approximated by
\begin{equation}
\int_{\left(\left]x^\pm_{i-1,j},\,x^\pm_{i,j}\right[\times \Delta^2_j\right)\cap \Gamma}
\phi^\pm(t,\mathbf x)\,n^{1\pm}(\mathbf x)\,\mathrm dS
\simeq
\bigl(A^{1\pm}_{i-\sfrac12,j}-B^{1\pm}_{i,j}\bigr)\,\Phi^{\gamma\pm}_{i,j}(t)
-\bigl(A^{1\pm}_{i-\sfrac12,j}-B^{1\pm}_{i-1,j}\bigr)\,\Phi^{\gamma\pm}_{i-1,j}(t),
\label{eq:interface-grad}
\end{equation}
leading to
\[
\operatorname{grad}^{1\gamma}_{i-\sfrac12,j}(A,B,W,\Phi)
:= \frac{
\bigl(A_{i-\sfrac12,j}-B_{i,j}\bigr)\Phi_{i,j}
-\bigl(A_{i-\sfrac12,j}-B_{i-1,j}\bigr)\Phi_{i-1,j}
}{W_{i-\sfrac12,j}}.
\]

The second component $G^{2\pm}_{i,j-\sfrac12}(t)$ is constructed analogously on
$W^{2\pm}_{i,j-\sfrac12}$ using $B^{2\pm}$ and $A^{2\pm}$:
\begin{align*}
\operatorname{grad}^{2\omega}_{i,j-\sfrac12}(B,W,\Phi)
&:= \frac{B_{i,j}\Phi_{i,j}-B_{i,j-1}\Phi_{i,j-1}}{W_{i,j-\sfrac12}},
\\
\operatorname{grad}^{2\gamma}_{i,j-\sfrac12}(A,B,W,\Phi)
&:= \frac{
\bigl(A_{i,j-\sfrac12}-B_{i,j}\bigr)\Phi_{i,j}
-\bigl(A_{i,j-\sfrac12}-B_{i,j-1}\bigr)\Phi_{i,j-1}
}{W_{i,j-\sfrac12}}.
\end{align*}

Finally, the constitutive relation Eq.~\eqref{eq:constitutive} is imposed in
discrete form by
\[
Q^{1\pm}_{i-\sfrac12,j}(t) \simeq -K^\pm\,G^{1\pm}_{i-\sfrac12,j}(t),
\qquad
Q^{2\pm}_{i,j-\sfrac12}(t) \simeq -K^\pm\,G^{2\pm}_{i,j-\sfrac12}(t),
\]
where $K^\pm$ is constant in each phase.

\subsection{Semi-discrete interface conditions}
\label{sec:semi_discrete_interface}

The interfacial conservation law
Eq.~\eqref{eq:surface-balance} corresponds to continuity of the
diffusive flux,
\[
\llbracket \mathbf q\cdot \mathbf n \rrbracket = 0 \quad \text{on }\Gamma.
\]
Integrating over each local interface segment $\Gamma_{i,j}$ yields
\begin{equation}
\int_{\Gamma_{i,j}} \llbracket \mathbf q\cdot \mathbf n \rrbracket \,\mathrm dS = 0.
\label{eq:semi_int_fluxbalance}
\end{equation}

The integrals in Eq.~\eqref{eq:semi_int_fluxbalance} correspond to the interfacial
contributions appearing in the Gauss decomposition in Eq.~\eqref{eq:gauss}. They are
approximated using the discrete interfacial divergence operators
$\operatorname{div}^{\alpha\gamma}$ introduced in Eq.~\eqref{eq:div-discrete-gamma}.
Consequently, the flux balance is semi-discretized as
\begin{equation}
\sum_{\alpha=1}^{d}
\Big[
V^+_{i,j}\,\operatorname{div}^{\alpha\gamma}_{i,j}\!\left(A^{\alpha+},B^{\alpha+},Q^{\alpha+}(t)\right)
-
V^-_{i,j}\,\operatorname{div}^{\alpha\gamma}_{i,j}\!\left(A^{\alpha-},B^{\alpha-},Q^{\alpha-}(t)\right)
\Big]
=0.
\label{eq:semi_int_fluxbalance_discrete}
\end{equation}

For mixed cells $(i,j)$, the jump/closure condition in Eq.~\eqref{eq:continuity}
(or in Eq.~\eqref{eq:henry}) is semi-discretized as
\begin{equation}
\Phi^{\gamma+}_{i,j}(t) - \lambda\,\Phi^{\gamma-}_{i,j}(t) = F_{i,j}(t),
\label{eq:semi_int_jump_discrete}
\end{equation}
where
\begin{equation*}
F_{i,j}(t)
:=
\frac{\langle \Gamma_{i,j},\, f(\cdot,t) \rangle}{\langle \Gamma_{i,j}, 1 \rangle}
=
\frac{\int_{\Gamma_{i,j}} f(\mathbf x,t)\,\mathrm dS}{\int_{\Gamma_{i,j}} \mathrm dS}.
\label{eq:Favg}
\end{equation*}
and
$\lambda$ is the interfacial coefficient appearing in the weighted jump
definition Eq.~\eqref{eq:weighted-jump} (e.g.\ $\lambda= 1$ adn $f=0$ to impose continuity of $\phi$ across the interface).

\subsection{Discrete equations}
\label{sec:discrete_equations}

We now use  a $\theta$-scheme time discretization of the semi-discrete bulk
balances Eq.~\eqref{eq:bulk-balance-semi}. For $0\le \theta \le 1$, we approximate
fluxes and sources at $t_{n+\theta}:=(1-\theta)t_n+\theta t_{n+1}$:
\begin{equation}
V^\pm_{i,j} C^\pm \frac{\Phi^{\omega\pm}_{n+1,i,j}-\Phi^{\omega\pm}_{n,i,j}}{\Delta t_n}
+
\sum_{\alpha=1}^{d}
V^\pm_{i,j}\,\operatorname{div}^{\alpha}_{i,j}\!\left(B^{\alpha\pm},Q^{\alpha\pm}_{n+\theta}\right)
=
V^\pm_{i,j}\,R^\pm_{n+\theta,i,j},
\label{eq:theta_bulk}
\end{equation}
where 
\begin{equation*}
R_{i,j}^{\pm}(t)
:=
\frac{\langle \Omega_{i,j}^{\pm},\, r^{\pm}(t,\cdot) \rangle}{V_{i,j}^{\pm}}
=
\frac{\int_{\Omega_{i,j}^{\pm}} r^{\pm}(t,\mathbf x)\,\mathrm{d}V}
{\int_{\Omega_{i,j}^{\pm}} \mathrm{d}V},
\label{eq:Ravg}
\end{equation*}
and $\Delta t_n := t_{n+1}-t_n$. The choices $\theta=0$
and $\theta=1$ correspond to the forward and backward Euler schemes respectively and
$\theta=\tfrac12$ to the midpoint rule.

The bulk equations are closed by specifying the face fluxes using the discrete
constitutive relation :
\begin{equation}
Q^{\alpha\pm}_{n+\theta}
=
-\,K^\pm\Big[
\operatorname{grad}^{\alpha\omega}\!\left(
B^{\alpha\pm},W^{\alpha\pm},(1-\theta)\Phi^{\omega\pm}_{n}+\theta\Phi^{\omega\pm}_{n+1}
\right)
+
\operatorname{grad}^{\alpha\gamma}\!\left(
A^{\alpha\pm},B^{\alpha\pm},W^{\alpha\pm},\Phi^{\gamma\pm}_{n+\theta}
\right)
\Big].
\label{eq:theta_fluxclosure}
\end{equation}

For mixed cells, the interfacial flux balance is discretized consistently at
$t_{n+\theta}$:
\begin{equation}
\sum_{\alpha=1}^{d}
\Big[
V^+_{i,j}\,\operatorname{div}^{\alpha\gamma}_{i,j}\!\left(A^{\alpha+},B^{\alpha+},Q^{\alpha+}_{n+\theta}\right)
-
V^-_{i,j}\,\operatorname{div}^{\alpha\gamma}_{i,j}\!\left(A^{\alpha-},B^{\alpha-},Q^{\alpha-}_{n+\theta}\right)
\Big]
=0.
\label{eq:theta_int_fluxbalance}
\end{equation}
Finally, the jump/closure condition is enforced at $t_{n+\theta}$ (where
$\Phi^{\gamma\pm}_{n+\theta,i,j}$ are collocated):
\begin{equation}
\Phi^{\gamma+}_{n+\theta,i,j} - \lambda\,\Phi^{\gamma-}_{n+\theta,i,j}
= F_{n+\theta,i,j}.
\label{eq:theta_int_jump}
\end{equation}

Choosing $\Phi^{\gamma\pm}$ at $t_{n+\theta}$ as primary unknowns (as opposed to
$t_n$) avoids the need to explicitly extrapolate interfacial values when
$\theta\neq 1$. Other time integration schemes may be employed, provided that
bulk and interfacial conservation statements are discretized consistently.

\subsection{Block structure of the linear system}
\label{sec:block_system}

Since the continuum model is linear, the corresponding discrete equations are linear and can
be written in matrix-vector form. Even for $\theta=0$ (forward Euler), the
unknowns at time level $n+1$ remain coupled through the interfacial conditions.
The matrices are nevertheless sparse, as the discrete operators were designed
to keep the stencil compact. The numerical solution at each step therefore
amounts to solving a sparse linear system coupling bulk and interfacial
unknowns.

We first eliminate the face fluxes $Q^{\alpha\pm}$ using the discrete
constitutive relation Eq.~\eqref{eq:theta_fluxclosure}, leaving
$\Phi^{\omega\pm}$ and $\Phi^{\gamma\pm}$ as the remaining unknowns. Let the
(diagonal) mass matrices be defined by
\[
M^\pm := \frac{C^\pm}{\Delta t_n}\,V^\pm,
\]
where $V^\pm$ denotes the diagonal matrix of phase volumes $V^\pm_{i,j}$.

Let $\operatorname{div}^{\alpha}$ denote the full (bulk+interface) directional
divergence operator used in Eq.~\eqref{eq:theta_bulk} and let
$\operatorname{grad}^{\alpha\omega}$ and $\operatorname{grad}^{\alpha\gamma}$
denote the bulk and interfacial contributions to the directional gradient.
Since these operators are linear in their arguments, we introduce the linear operators
\begin{align}
L^{\omega\omega,\pm}
&:= -\sum_{\alpha\in\{1,2\}}
\partial_{Q}\bigl(V^\pm\,\operatorname{div}^{\alpha}(B^{\alpha\pm},\,\cdot)\bigr)\;
K^{\alpha\pm}\;
\partial_{\Phi}\operatorname{grad}^{\alpha\omega}(B^{\alpha\pm},W^{\alpha\pm},\,\cdot),
\label{eq:Lww}\\
L^{\omega\gamma,\pm}
&:= -\sum_{\alpha\in\{1,2\}}
\partial_{Q}\bigl(V^\pm\,\operatorname{div}^{\alpha}(B^{\alpha\pm},\,\cdot)\bigr)\;
K^{\alpha\pm}\;
\partial_{\Phi}\operatorname{grad}^{\alpha\gamma}(A^{\alpha\pm},B^{\alpha\pm},W^{\alpha\pm},\,\cdot),
\label{eq:Lwg}\\
L^{\gamma\omega,\pm}
&:= -\sum_{\alpha\in\{1,2\}}
\partial_{Q}\bigl(V^\pm\,\operatorname{div}^{\alpha\gamma}(A^{\alpha\pm},B^{\alpha\pm},\,\cdot)\bigr)\;
K^{\alpha\pm}\;
\partial_{\Phi}\operatorname{grad}^{\alpha\omega}(B^{\alpha\pm},W^{\alpha\pm},\,\cdot),
\label{eq:Lgw}\\
L^{\gamma\gamma,\pm}
&:= -\sum_{\alpha\in\{1,2\}}
\partial_{Q}\bigl(V^\pm\,\operatorname{div}^{\alpha\gamma}(A^{\alpha\pm},B^{\alpha\pm},\,\cdot)\bigr)\;
K^{\alpha\pm}\;
\partial_{\Phi}\operatorname{grad}^{\alpha\gamma}(A^{\alpha\pm},B^{\alpha\pm},W^{\alpha\pm},\,\cdot).
\label{eq:Lgg}
\end{align}
Equivalently, one may view these as Jacobians; here they are constant matrices.

The discrete system at a given time step can then be written as
\begin{multline}
\left[
\begin{array}{c|c|cc}
M^- + \theta L^{\omega\omega,-} & 0 &  L^{\omega\gamma,-} & 0 \\ \hline
0 & M^+ + \theta L^{\omega\omega,+} & 0 &  L^{\omega\gamma,+} \\ \hline
-\theta L^{\gamma\omega,-} & \theta L^{\gamma\omega,+} & - L^{\gamma\gamma,-} & L^{\gamma\gamma,+} \\
0 & 0 & -\lambda I & I
\end{array}
\right]
\left[
\begin{array}{c}
\Phi^{\omega-}_{n+1} \\
\Phi^{\omega+}_{n+1} \\
\Phi^{\gamma-}_{n+\theta} \\
\Phi^{\gamma+}_{n+\theta}
\end{array}
\right] \\
=
\left[
\begin{array}{cc}
M^- - (1-\theta)L^{\omega\omega,-} & 0 \\
0 & M^+ - (1-\theta)L^{\omega\omega,+} \\
(1-\theta)L^{\gamma\omega,-} & -(1-\theta)L^{\gamma\omega,+} \\
0 & 0
\end{array}
\right]
\left[
\begin{array}{c}
\Phi^{\omega-}_{n} \\
\Phi^{\omega+}_{n}
\end{array}
\right]
+
\left[
\begin{array}{c}
V^- R^-_{n+\theta} \\
V^+ R^+_{n+\theta} \\
0 \\
F_{n+\theta}
\end{array}
\right].
\label{eq:block_system}
\end{multline}

The horizontal and vertical separators highlight the arrow-type structure of
the matrix, characteristic of domain-decomposition couplings between subdomain
unknowns and interfacial unknowns.

Finally, for symmetric positive definite face-centered diffusion coefficients
$K^{\alpha\pm}$, the operators $L^{\omega\omega,\pm}$ are symmetric positive definite on the active bulk unknowns. Moreover, under the standard
finite-volume duality between divergence and gradient, the
discrete adjoint relation also holds, 
\[
-\,L^{\gamma\omega,\pm} = \left(L^{\omega\gamma,\pm}\right)^{\top}.
\]
The resulting sparse linear system Eq.~\eqref{eq:block_system} is solved using an
appropriate linear solver (direct or iterative), depending on problem size.
A one-dimensional example is presented in App.~\ref{appendix:1D_diffusion}.

\section{Numerical validation in static domains}
\label{sec:numstatic}

We validate the accuracy and robustness of the proposed cut-cell formulation
for static embedded geometries. Both single-phase and two-phase
configurations are considered, for elliptic (steady diffusion)
and parabolic (unsteady diffusion) problems. Analytical solutions and
corresponding forcing terms used for verification are provided in each paragraph.

Let $\phi^{\rm ex}(\mathbf x,t)$ denote the exact solution. Errors are evaluated from
cell-averaged unknowns at selected times $t_n$ (typically the final time),
and reported separately over regular cells and cut cells.
The index sets are  defined by
\[
\mathcal I_{\rm reg}=\{(i,j):\ \Gamma_{i,j}=\emptyset\},\qquad
\mathcal I_{\rm cut}=\{(i,j):\ \Gamma_{i,j}\neq\emptyset\},\qquad
\mathcal I_{\rm all}=\mathcal I_{\rm reg}\cup \mathcal I_{\rm cut}.
\]

For any subset $S\in\{\mathrm{reg},\mathrm{cut},\mathrm{all}\}$, we define the
discrete $L^2$ error norm as
\begin{equation}
\|e(t_n)\|_{2,S}^{\rm}
=
\left(
\frac{\displaystyle
\sum_{(i,j)\in\mathcal I_S}
V_{i,j}\,
\bigl|\Phi^{\rm \omega}_{n,i,j} - \Phi^{\rm ex}_{n,i,j}\bigr|^2}
{\displaystyle
\sum_{(i,j)\in\mathcal I_S}
V_{i,j}\,
}
\right)^{1/2},
\label{eq:relL2}
\end{equation}
where 
$\phi^{\rm ex}_{n,i,j}$ denotes the exact solution evaluated consistently with
the discrete unknown e.g.\ at the cell centroid for cell averages. We report
\[
\|e\|_{2,\mathrm{reg}}^{\rm},\qquad
\|e\|_{2,\mathrm{cut}}^{\rm},\qquad
\|e\|_{2,\mathrm{all}}^{\rm}.
\]

The empirical convergence order is estimated between two successive grid
resolutions $h_i$ and $h_{i+1}$ using
\begin{equation}
p_{2,S} =
\frac{\log\!\left(\|e\|_{2,S}^{(i)} / \|e\|_{2,S}^{(i+1)}\right)}
     {\log\!\left(h_i / h_{i+1}\right)},
\qquad
S \in \{\mathrm{reg},\mathrm{cut},\mathrm{all}\}.
\label{eq:order}
\end{equation}

In addition, a global convergence rate is computed as the least-squares
slope of $\log(\|e\|_{2,S})$ versus $\log(h)$ over all considered grid
resolutions.

\subsection{Single-phase validation}
\label{sec:mono_validation}

We first assess the single-phase formulation in a static embedded geometry.
The single-phase formulation is obtained by restricting the computation to the
light phase $\Omega^-$, i.e.\ by solving only for $\Phi^{\omega-}$ and
$\Phi^{\gamma-}$ and dropping all ``$+$'' variables and interfacial coupling
equations. In the continuous setting, the embedded boundary
$\Gamma=\partial\Omega^-\cap\Omega$ is equipped with the Robin condition
\[
\mathbf q^-\!\cdot\mathbf n + \beta\,\phi^- = f
\qquad \text{on }\Gamma,
\]
where $\mathbf n$ is the outward unit normal to $\Omega^-$. Integrating over a
boundary segment $\Gamma_{i,j}$ gives
\[
\int_{\Gamma_{i,j}} \mathbf q^-\!\cdot\mathbf n\,\mathrm dS
+\beta\int_{\Gamma_{i,j}}\phi^-\,\mathrm dS
=
\int_{\Gamma_{i,j}} f\,\mathrm dS.
\]
Using the interfacial (cut) operator $\operatorname{div}^{\alpha\gamma}$ to
represent the boundary-flux contribution consistently with the bulk balance,
the discrete Robin condition is enforced as
\[
\sum_{\alpha=1}^{d} V^-_{i,j}\,
\operatorname{div}^{\alpha\gamma}_{i,j}\!\left(A^{\alpha-},B^{\alpha-},Q^{\alpha-}_{n+\theta}\right)
\;+\;
\beta\,\langle \Gamma_{i,j},1\rangle\,\Phi^{\gamma-}_{n+\theta,i,j}
\;=\;
\langle \Gamma_{i,j},1\rangle\,F_{n+\theta,i,j},
\]
with the segment-averaged value $F_{n+\theta,i,j}$.

Dirichlet or Neumann boundary conditions are recovered in the limits $\beta\to\infty$ (with $f=\beta\,\phi_D$) and $\beta=0$ respectively.

\subsubsection{Johansen-Colella Problem~1: Star-shaped Poisson problem}
\label{sec:jc_p1}

We reproduce the classical embedded-boundary verification problem proposed by
Johansen and Colella~\cite{johansen_cartesian_1998}. The physical domain
$\Omega$ is star-shaped and defined in polar coordinates by
\begin{equation}
\Omega=\{(r,\theta):\ 0\le r \le R(\theta)\},
\qquad
R(\theta)=0.30+0.15\cos(6\theta).
\label{eq:jc_domain}
\end{equation}
We solve the Poisson equation
\begin{equation}
\Delta \phi = 7r^2\cos(3\theta)
\qquad \text{in }\Omega,
\label{eq:jc_poisson}
\end{equation}
supplemented with Dirichlet boundary conditions taken from the exact solution
\begin{equation}
\phi(r,\theta)=r^4\cos(3\theta)
\qquad \text{on }\Gamma.
\label{eq:jc_exact}
\end{equation}
The right-hand side in Eq.~\eqref{eq:jc_poisson} follows from the polar Laplacian
identity $\Delta(r^m\cos(k\theta))=(m^2-k^2)\,r^{m-2}\cos(k\theta)$ with
$m=4$ and $k=3$.

The problem is discretized on a uniform Cartesian grid of spacing $h$ using the cut-cell operators of Sec.~\ref{sec:cutcell}, with $\Gamma$ embedded
in the background mesh. Errors are reported using the $L^2$ norms
defined in Eq.~\eqref{eq:relL2}, distinguishing regular cells (fully inside
$\Omega$), cut cells (intersected by $\Gamma$) and all active cells.

\begin{figure}[!h]
  \centering
  \begin{subfigure}[b]{0.49\linewidth}
    \centering
    \includegraphics[width=\linewidth]{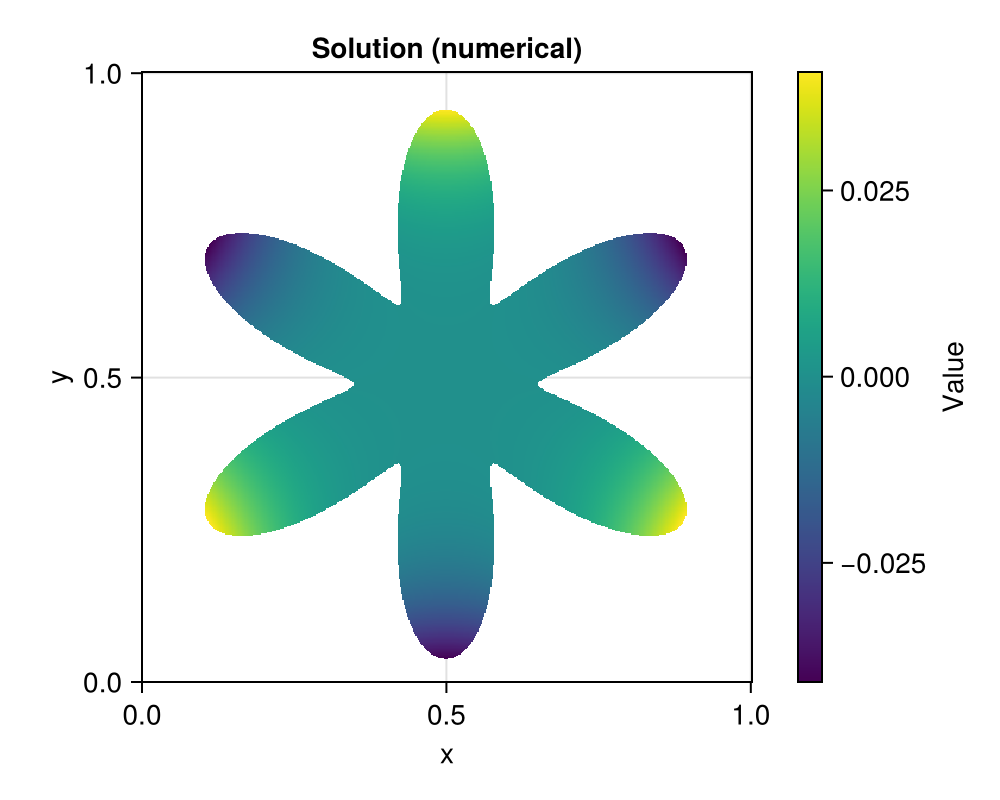}
    \caption{Numerical solution $\phi^{\rm num}$}
    \label{fig:jc_p1_sol}
  \end{subfigure}
  \hfill
  \begin{subfigure}[b]{0.49\linewidth}
    \centering
    \includegraphics[width=\linewidth]{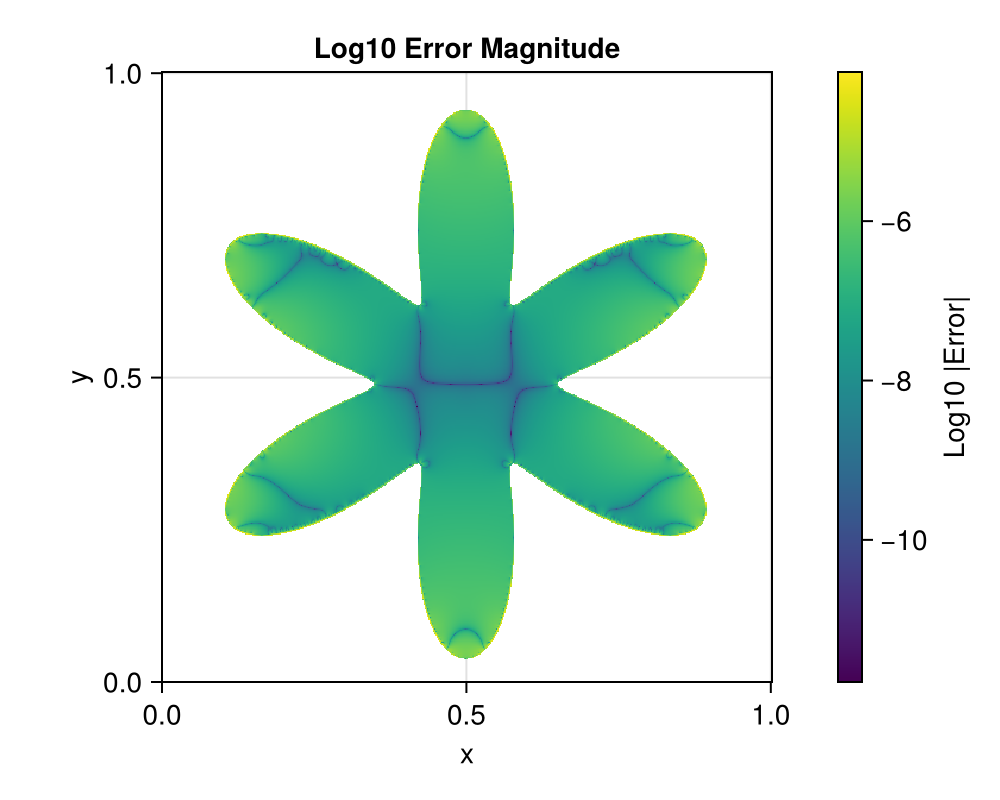}
    \caption{Logarithm of the pointwise error $\log_{10} |\phi^{\rm num}-\phi^{\rm ex}|$}
    \label{fig:jc_p1_err}
  \end{subfigure}
  \caption{Johansen-Colella Problem~1 on the star-shaped domain
  Eq.~\eqref{eq:jc_domain}: numerical solution and corresponding $\log_{10}$ error field on a $512^2$ mesh}
  \label{fig:jc_p1}
\end{figure}

\begin{table}[!h]
\centering
\setlength{\tabcolsep}{6pt}
\renewcommand{\arraystretch}{1.15}
\begin{tabular}{c c c | c c c | c c c}
\hline
$h$ & $N_x{=}N_y$ & $N_{\Omega}$ &
$\|e\|_{2,\mathrm{all}}$ &
$\|e\|_{2,\mathrm{reg}}$ &
$\|e\|_{2,\mathrm{cut}}$ &
$p_{\mathrm{all}}$ & $p_{\mathrm{reg}}$ & $p_{\mathrm{cut}}$ \\
\hline
0.0625      & 15  & 39    & $3.726e{-3}$ & $1.236e{-3}$ & $3.515e{-3}$ & -   & -   & -   \\
0.03125     & 29  & 228   & $1.498e{-3}$ & $1.067e{-3}$ & $1.052e{-3}$ & 1.31 & 0.21 & 1.74 \\
0.015625    & 58  & 1138  & $3.098e{-5}$ & $1.504e{-5}$ & $2.708e{-5}$ & 5.60 & 6.15 & 5.28 \\
0.0078125   & 116 & 4868  & $1.139e{-5}$ & $4.821e{-6}$ & $1.032e{-5}$ & 1.44 & 1.64 & 1.39 \\
0.00390625  & 231 & 20166 & $3.958e{-6}$ & $1.399e{-6}$ & $3.703e{-6}$ & 1.52 & 1.79 & 1.48 \\
0.001953125 & 461 & 82018 & $1.327e{-6}$ & $4.054e{-7}$ & $1.264e{-6}$ & 1.58 & 1.79 & 1.55 \\
\hline
fit   & -  & - & - & - & - &  2.41 & 2.52 & 2.37 \\
\hline
\end{tabular}
\caption{Johansen-Colella Problem~1: $L^2$ errors in regular, cut and
all active cells and pairwise convergence orders. Here $N_\Omega$ denotes the
number of active cells intersecting $\Omega$ (i.e.\ $V^-_{i,j}>0$).}
\label{tab:jc_p1}
\end{table}

On the finest grids, the observed convergence rates stabilize around
$p_{\mathrm{all}}\approx 1.6$, with slightly higher rates in regular cells
$p_{\mathrm{reg}}\approx 1.8$ and comparable rates in cut cells
$p_{\mathrm{cut}}\approx 1.5$. The unusually large pairwise orders reported
between some intermediate resolutions (e.g.\ $p\gtrsim 5$) reflect a clearly
pre-asymptotic regime, in which the star-shaped boundary is still
under-resolved (small $N_\Omega$) and the error is dominated by geometry-related
effects and occasional favorable cancellations rather than the asymptotic truncation
error.

This trend is consistent with the error fields shown in
Fig.~\ref{fig:jc_p1}: the point-wise error concentrates near the embedded
boundary and is larger in cut cells, while the interior region exhibits
a smoother and smaller error distribution. As the grid is refined, the boundary
layer of error contracts and the regular-cell contribution approaches the nominal
second-order behavior, whereas the global rate remains slightly reduced due to
the persistent influence of the cut-cell closure at the curved interface.

\subsubsection{Johansen-Colella Problem 2: Flower-shaped steady diffusion (boundedness test).}

As an additional boundedness/maximum-principle diagnostic, we reproduce a
steady diffusion problem in the spirit of Johansen and Colella~\cite{johansen_cartesian_1998}. We consider a steady
diffusion problem in the unit box with an embedded flower-shaped hole.
The computational domain is the fixed Cartesian box
\(
\Omega_{\rm box} = [0,1]\times[0,1],
\)
from which a flower-shaped region is removed. The immersed boundary $\Gamma$ is
defined in polar coordinates (centered at $\mathbf x_c=(0.5,0.5)$) by
\begin{equation}
r = R(\theta) = 0.25 + 0.05\cos(6\theta),
\label{eq:jc_flower}
\end{equation}
and the physical domain is the exterior region
\[
\Omega = \Omega_{\rm box}\setminus \{(r,\theta): 0\le r \le R(\theta)\}.
\]
We solve the Laplace equation
\begin{equation}
\Delta \phi = 0
\qquad \text{in }\Omega,
\label{eq:jc_flower_laplace}
\end{equation}
with Dirichlet boundary conditions prescribed on both boundaries:
\begin{equation}
\phi = 1 \quad \text{on }\Gamma,
\qquad
\phi = 0 \quad \text{on }\partial\Omega_{\rm box}.
\label{eq:jc_flower_bc}
\end{equation}
By the maximum principle, the exact solution is bounded within $[0,1]$. We
therefore monitor the discrete extrema $\min\phi^h$ and $\max\phi^h$ and the
overshoot/undershoot ratio $k/N$, where $k$ is the number of active cells with
values outside $[0,1]$ and $N$ is the total number of active cells. No overshoot
or undershoot is observed at any resolution, as summarized in
Tab.~\ref{tab:johansen-flower}. To stress the robustness of the proposed method, extremely coarse resolution are also considered here (only $4\times4$
cells).

\begin{figure}[!h]
\centering
\includegraphics[width=0.5\linewidth]{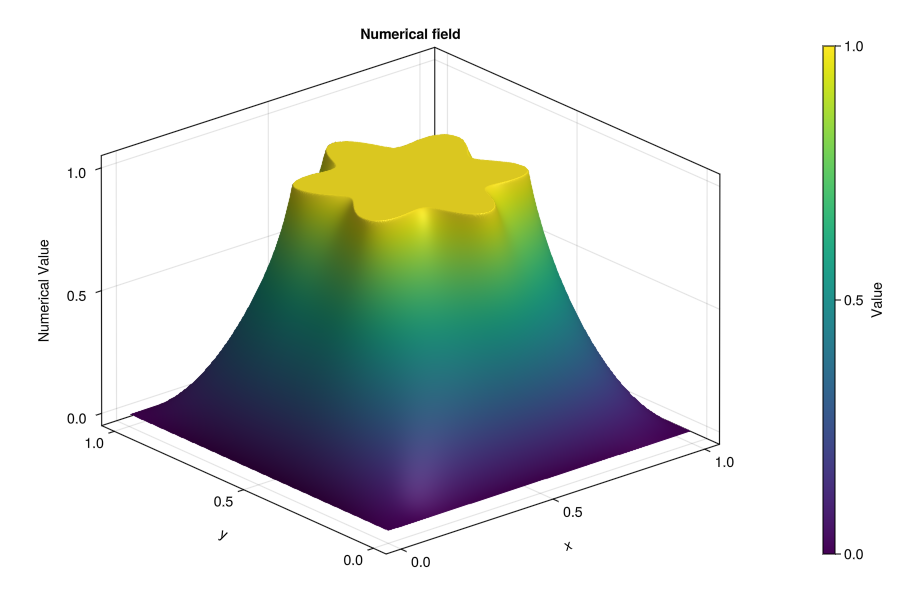}
\caption{Johansen-Colella Problem~2: numerical steady diffusion field
$\phi^{\rm num}$ in the unit box, with the embedded flower-shaped boundary
Eq.~\eqref{eq:jc_flower} defining a hole.}
\label{fig:jc_p2_field}
\end{figure}

\begin{table}[h!]
\centering
\setlength{\tabcolsep}{6pt}
\renewcommand{\arraystretch}{1.2}
\begin{tabular}{c c | c c | c c c}
\hline
$N_x$ & $N_y$ &
$N_{\mathrm{cut}}$ & $N_{\mathrm{inside}}$ &
$\max \phi^h$ & $\min \phi^h$ &
$k/N$ \\
\hline
4   & 4   & 8   & 7     & 1.0 & 0.0 & 0.0 \\
8   & 8   & 16  & 43    & 1.0 & 0.0 & 0.0 \\
16  & 16  & 40  & 183   & 1.0 & 0.0 & 0.0 \\
32  & 32  & 90  & 773   & 1.0 & 0.0 & 0.0 \\
64  & 64  & 160 & 3195  & 1.0 & 0.0 & 0.0 \\
128 & 128 & 328 & 12939 & 1.0 & 0.0 & 0.0 \\
256 & 256 & 664 & 52099 & 1.0 & 0.0 & 0.0 \\
\hline
\end{tabular}
\caption{Johansen-Colella flower-shaped steady diffusion (boundedness test):
discrete extrema and overshoot ratio $k/N$ for increasing grid resolution. No
violations of the physical bounds $[0,1]$ are detected.}
\label{tab:johansen-flower}
\end{table}

\subsubsection{Poisson equation in a disk with Robin boundary condition.}

We next validate the single-phase formulation for an elliptic problem with a
Robin boundary condition imposed on a curved embedded boundary. The
computational box is
\(
\Omega_{\rm box}=[0,4]\times[0,4],
\)
and the physical domain is the disk
\[
\Omega=\{(x,y)\in\Omega_{\rm box}:\ (x-x_c)^2+(y-y_c)^2< R^2\},
\qquad
(x_c,y_c)=(2,2),\quad R=1.
\]
We solve the Poisson equation
\begin{equation}
-\Delta \phi = 1
\qquad \text{in }\Omega,
\label{eq:robin_disk_poisson}
\end{equation}
supplemented with a Robin boundary condition on the embedded boundary
$\Gamma=\partial\Omega$,
\begin{equation}
\partial_n \phi + \phi = 1
\qquad \text{on }\Gamma,
\label{eq:robin_disk_bc}
\end{equation}
where $\partial_n \phi=\nabla \phi\cdot\mathbf n$ and $\mathbf n$ denotes the outward
unit normal to $\Omega$. The exact solution is
\begin{equation}
\phi^{\rm ex}(x,y)=\frac{7}{4}-\frac{(x-x_c)^2+(y-y_c)^2}{4},
\label{eq:robin_disk_exact}
\end{equation}
which satisfies $-\Delta \phi^{\rm ex}=1$ in $\Omega$ and the Robin condition
Eq.~\eqref{eq:robin_disk_bc} on $\Gamma$.

\begin{figure}[!h]
  \centering
  \begin{subfigure}[b]{0.49\linewidth}
    \centering
    \includegraphics[width=\linewidth]{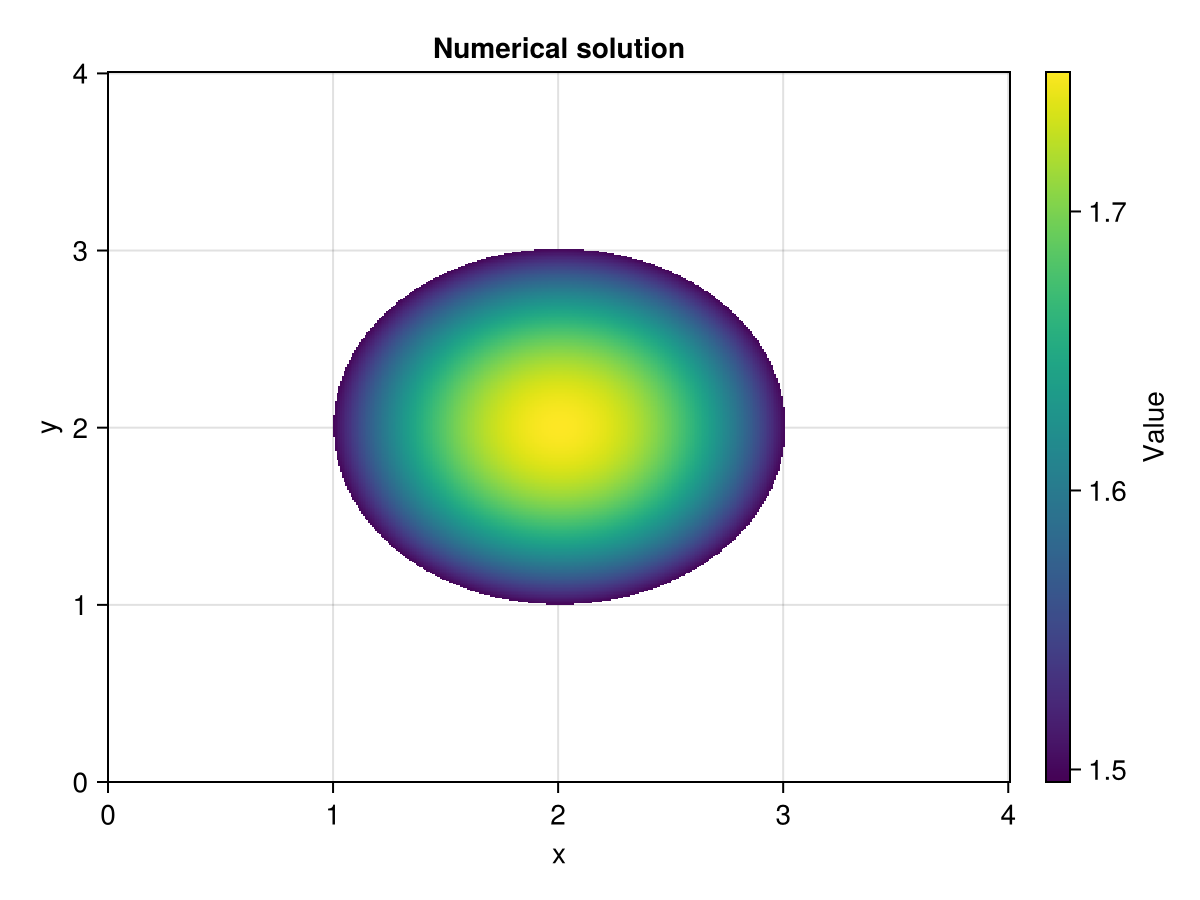}
    \caption{Numerical solution $u^{\rm num}$}
    \label{fig:robin_disk_sol}
  \end{subfigure}
  \hfill
  \begin{subfigure}[b]{0.49\linewidth}
    \centering
    \includegraphics[width=\linewidth]{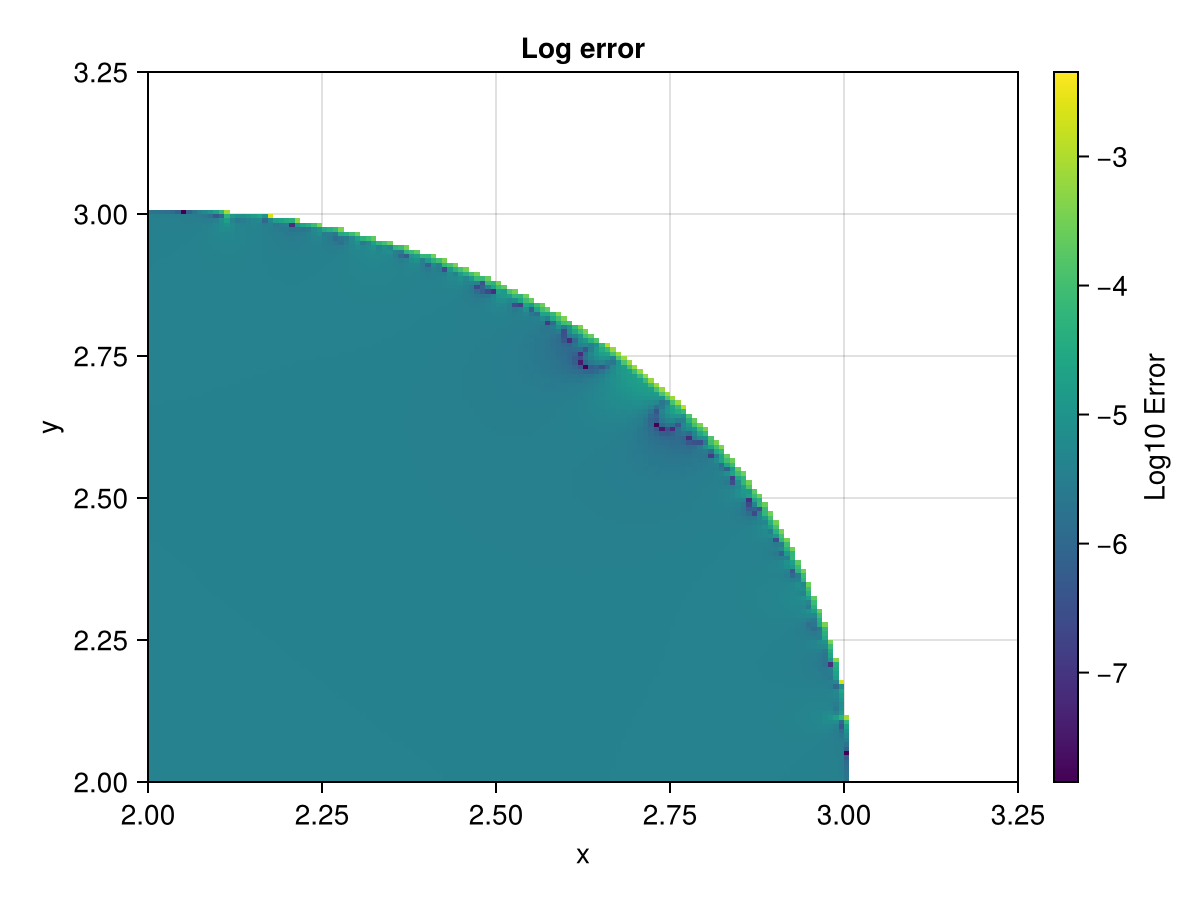}
    \caption{$\log_{10}\!\left(|u^{\rm num}-u^{\rm ex}|\right)$}
    \label{fig:robin_disk_logerr}
  \end{subfigure}
  \caption{Poisson problem in a disk with Robin boundary condition
  Eq.~\eqref{eq:robin_disk_bc}: numerical solution and base-10 logarithm of the
  pointwise error.}
  \label{fig:robin_disk}
\end{figure}

 The discrete semi-norm $H^1$ error is defined from the gradient of the
numerical solution,
\[
\|e\|_{H^1,{S}}^{\rm}
=
\left(
\frac{\displaystyle\sum_{(i,j)\in\mathcal I_{S}}
V_{i,j}\,\bigl\|\nabla \Phi_{i,j}-\nabla \Phi^{\rm ex}_{i,j}\bigr\|^2}
{\displaystyle\sum_{(i,j)\in\mathcal I_{S}}
V_{i,j}\,}
\right)^{1/2},
\]
where $\nabla \Phi$ is obtained with the discrete gradient operator of
Sec.~\ref{sec:discrete_bulk} and $\nabla \Phi^{\rm ex}$ is evaluated at consistent
locations. Pairwise orders are computed using Eq.~\eqref{eq:order}. Tab.~\ref{tab:robin_disk_conv} reports the $L^2$ errors evaluated over all, regular and cut cells. Tab.~\ref{tab:robin_disk_H1} reports the $H^1$ semi-norm errors evaluated over
all, regular and cut cells.

\begin{table}[!h]
\centering
\setlength{\tabcolsep}{6pt}
\renewcommand{\arraystretch}{1.15}
\begin{tabular}{c c | c c c | c c c}
\hline
$h$ & $N_x{=}N_y$ &
$\|e\|_{2,\mathrm{all}}$ &
$\|e\|_{2,\mathrm{reg}}$ &
$\|e\|_{2,\mathrm{cut}}$ &
$p_{2,\mathrm{all}}$ &
$p_{2,\mathrm{reg}}$ &
$p_{2,\mathrm{cut}}$ \\
\hline
0.25      & 8   & $5.413e{-3}$ & $3.883e{-3}$ & $3.771e{-3}$ & -    & -    & -    \\
0.125     & 16  & $1.216e{-3}$ & $8.813e{-4}$ & $8.382e{-4}$ & 2.15 & 2.14 & 2.17 \\
0.0625    & 32  & $4.198e{-4}$ & $2.616e{-4}$ & $3.284e{-4}$ & 1.53 & 1.75 & 1.35 \\
0.03125   & 64  & $1.209e{-4}$ & $7.104e{-5}$ & $9.788e{-5}$ & 1.80 & 1.88 & 1.75 \\
0.015625  & 128 & $3.476e{-5}$ & $1.697e{-5}$ & $3.033e{-5}$ & 1.80 & 2.07 & 1.69 \\
0.0078125 & 256 & $1.000e{-5}$ & $4.056e{-6}$ & $9.142e{-6}$ & 1.80 & 2.06 & 1.73 \\
\hline
fit & - & - & - & - & 1.79 & 1.96 & 1.70 \\
\hline
\end{tabular}
\caption{$L^2$ error convergence for the Poisson problem in a disk with Robin boundary condition, measured over all, regular, and cut cells. }
\label{tab:robin_disk_conv}
\end{table}

\begin{table}[!h]
\centering
\setlength{\tabcolsep}{6pt}
\renewcommand{\arraystretch}{1.15}
\begin{tabular}{c c | c c c | c c c}
\hline
$h$ & $N_x{=}N_y$ &
$\|e\|_{H^1,\mathrm{all}}$ &
$\|e\|_{H^1,\mathrm{reg}}$ &
$\|e\|_{H^1,\mathrm{cut}}$ &
$p_{H^1,\mathrm{all}}$ &
$p_{H^1,\mathrm{reg}}$ &
$p_{H^1,\mathrm{cut}}$ \\
\hline
0.25      & 8   & $8.888e{-2}$ & $7.618e{-2}$ & $4.578e{-2}$ & -    & -    & -    \\
0.125     & 16  & $4.456e{-2}$ & $4.103e{-2}$ & $1.738e{-2}$ & 1.00 & 0.89 & 1.40 \\
0.0625    & 32  & $2.290e{-2}$ & $2.130e{-2}$ & $8.419e{-3}$ & 0.96 & 0.95 & 1.05 \\
0.03125   & 64  & $1.181e{-2}$ & $1.089e{-2}$ & $4.550e{-3}$ & 0.96 & 0.97 & 0.89 \\
0.015625  & 128 & $5.797e{-3}$ & $5.487e{-3}$ & $1.871e{-3}$ & 1.03 & 0.99 & 1.28 \\
0.0078125 & 256 & $2.909e{-3}$ & $2.763e{-3}$ & $9.086e{-4}$ & 0.99 & 0.99 & 1.04 \\
\hline
fit & - & - & - & - & 0.98 & 0.96 & 1.11 \\
\hline
\end{tabular}
\caption{$H^1$-semi-norm error convergence for the Poisson problem in a disk with Robin boundary condition, measured over all, regular, and cut cells.}
\label{tab:robin_disk_H1}
\end{table}

Tab.~\ref{tab:robin_disk_conv} confirms that the cut-cell formulation enforces
the Robin condition on the curved embedded boundary with the expected accuracy.
The $L^2$ error over regular cells converges at essentially
second order once the mesh is sufficiently fine ($p_{2,\mathrm{reg}}\approx 2$
for the last refinements), indicating that the scheme retains its nominal
accuracy away from the boundary. Over all cells, the observed rate
stabilizes around $p_{2,\mathrm{all}}\approx 1.8$, reflecting the fact that the
global error is influenced by the cut-cell region where geometric trimming and
boundary closure dominate. The $L^2$ error restricted to cut cells
converges slightly below second order ($p_{2,\mathrm{cut}}\approx 1.7$ on the
finest grids), which is consistent with the reduced regularity of the discrete
operators near embedded boundaries. Finally, Tab.~\ref{tab:robin_disk_H1} shows that the discrete $H^1$-semi-norm errors converge at first order, with overall slopes close to one
($p_{H^1,\mathrm{reg}}\approx 1$, $p_{H^1,\mathrm{all}}\approx 1$ and
$p_{H^1,\mathrm{cut}}\approx 1$), as expected for gradient reconstruction in the
vicinity of irregular control volumes and curved boundary closures.

\subsubsection{Unsteady diffusion in a sphere with Robin boundary condition (3D).}
\label{sec:3dunsteadrobin}

We finally consider a three-dimensional single-phase parabolic test in a spherical
domain with a Robin boundary condition. The physical domain is the ball
\(
\Omega=\{\mathbf x\in\mathbb R^3:\ \|\mathbf x\|<R\},
\)
embedded in a fixed Cartesian box. Let $r=\|\mathbf x\|$. We solve
\begin{equation}
\partial_t \phi = a\,\Delta\phi
\qquad \text{in }\Omega,
\label{eq:ball_diffusion}
\end{equation}
with the radial Laplacian
\(
\Delta \phi =
\frac{1}{r^2}\frac{\partial}{\partial r}
\left(
r^2\frac{\partial \phi}{\partial r}
\right),
\)
and the Robin condition on $\Gamma=\partial\Omega$,
\begin{equation}
\left.\frac{\partial \phi}{\partial r}\right|_{r=R}
+
k\,\phi(R,t) = 0.
\label{eq:ball_robin}
\end{equation}
Starting from a uniform initial condition $\phi(\mathbf x,0)=\phi_0$, an exact
radial series solution is available \cite{polyanin_handbook_2001}. The eigenvalues
$\{\mu_n\}_{n\ge 1}$ satisfy
\begin{equation}
\mu_n \cot\mu_n + kR - 1 = 0,
\label{eq:robin_eigs}
\end{equation}
and the solution reads
\begin{equation}
\phi(r,t)
=
\sum_{n=1}^{\infty}
C_n\,
\frac{\sin(\mu_n r / R)}{r}\,
\exp\!\left(- a \mu_n^2 t/R^2\right),
\label{eq:ball_series}
\end{equation}
with coefficients
\begin{equation}
C_n
=
\frac{2 k R^2 \phi_0}{\mu_n^2}
\,
\frac{\mu_n^2 + (kR - 1)^2}{\mu_n^2 + kR(kR - 1)}
\,
\sin(\mu_n).
\label{eq:ball_coeffs}
\end{equation}

Time integration is performed with the midpoint rule ($\theta=\tfrac12$) up to
$t_f=0.1$, using the time step
\(
\Delta t = 0.25\,\min(\Delta x^2,\Delta y^2,\Delta z^2).
\)
We report $L^2$ errors at $t_f$ over regular cells, cut cells and all
active cells. For convenience, we also indicate $N_{\rm diam}$, the number of
grid points across the sphere diameter, to highlight that the method remains
robust even on very coarse resolutions (down to $N_{\rm diam}=1$).

\begin{table}[h!]
\centering
\setlength{\tabcolsep}{6pt}
\renewcommand{\arraystretch}{1.2}
\begin{tabular}{c c | c c c | c c c}
\hline
$h$ & $N_{\rm diam}$ &
$\|e(t_f)\|_{2,\mathrm{reg}}$ &
$\|e(t_f)\|_{2,\mathrm{cut}}$ &
$\|e(t_f)\|_{2,\mathrm{all}}$ &
$p_{\mathrm{reg}}$ & $p_{\mathrm{cut}}$ & $p_{\mathrm{all}}$ \\
\hline
1.0     & 1  & $3.149e{-1}$ & $1.996e{-1}$ & $3.729e{-1}$ & -   & -   & -   \\
0.5     & 3  & $7.141e{-2}$ & $1.362e{-1}$ & $1.538e{-1}$ & 2.14 & 0.55 & 1.28 \\
0.25    & 9  & $3.582e{-2}$ & $3.316e{-2}$ & $4.881e{-2}$ & 1.00 & 2.04 & 1.66 \\
0.125   & 14 & $9.452e{-3}$ & $5.780e{-3}$ & $1.108e{-2}$ & 1.92 & 2.52 & 2.14 \\
0.0625  & 34 & $2.900e{-3}$ & $1.303e{-3}$ & $3.179e{-3}$ & 1.70 & 2.15 & 1.80 \\
0.03125 & 70 & $8.911e{-4}$ & $3.201e{-4}$ & $9.123e{-4}$ & 1.70 & 2.02 & 1.80 \\
\hline
fit     & - & - & - & - & 1.82 & 2.07 & 1.85 \\
\hline
\end{tabular}
\caption{$L^2$ errors at final time $t_f=0.1$ for the 3D single-phase
unsteady diffusion problem in a sphere with Robin boundary condition.}
\label{tab:3D-robin-final}
\end{table}

Tab.~\ref{tab:3D-robin-final} demonstrates that the proposed cut-cell framework
extends to three dimensions: the geometric moments (trimmed volumes,
face apertures and boundary measures) can be integrated robustly in 3D and
lead to stable discretizations even on very coarse grids (e.g.\ $N_{\rm diam}=1$).
On the finest refinements, the global error converges at approximately
second order, with fitted rates $p_{\mathrm{all}}\approx 1.9$ and
$p_{\mathrm{reg}}\approx 1.8$. The cut-cell error also exhibits near-second-order
behavior ($p_{\mathrm{cut}}\approx 2.0$ ), indicating that the Robin
closure and the moment-based representation of the curved boundary remain
accurate in three dimensions.

The non-monotone pairwise orders at coarse-to-intermediate resolutions
(e.g.\ the low $p_{\mathrm{cut}}$ between $h=1$ and $h=0.5$, followed by larger
rates) are symptomatic of a pre-asymptotic regime where the sphere is strongly
under-resolved and the error is dominated by boundary
closure effects. Once the diameter is resolved by a few tens of cells
($N_{\rm diam}\gtrsim 30$), the rates stabilize and the expected asymptotic
convergence is recovered.

\begin{figure}[!h]
  \centering
  \includegraphics[width=0.82\linewidth]{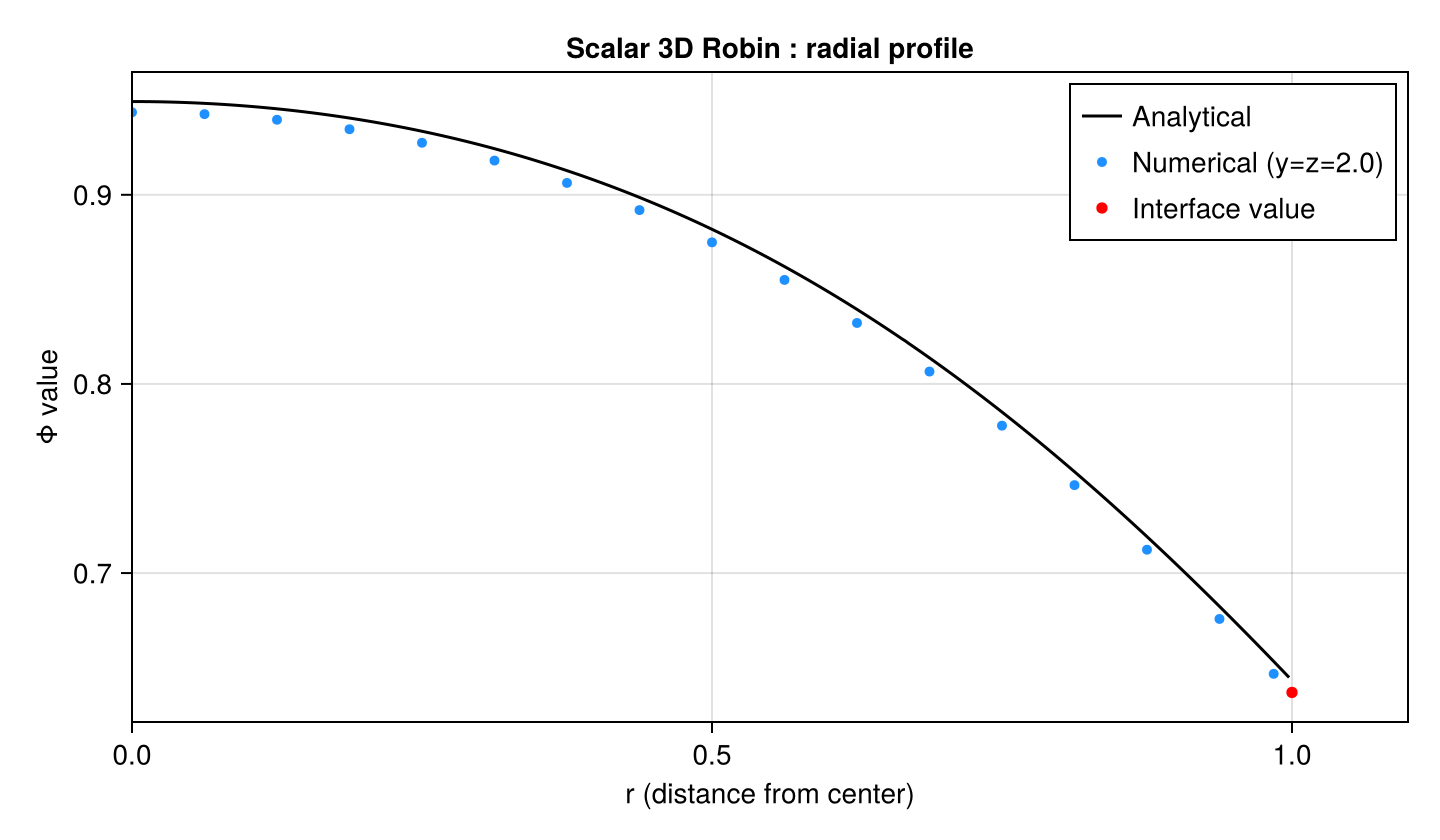}
  \caption{3D unsteady diffusion in a sphere with Robin boundary condition:
  radial profile of $\phi(r,t_f)$ along a line through the center.
  Blue markers: numerical samples as a function of $r$.
  Red marker: discrete boundary/interfacial value at $r=R$ used in the Robin
  closure. Solid line: analytical series solution Eq.~\eqref{eq:ball_series}.}
  \label{fig:3d_robin_radial_profile}
\end{figure}

Fig.~\ref{fig:3d_robin_radial_profile} shows the radial profile of $\phi(\cdot,t_f)$ extracted along a line passing
through the center. Blue markers denote the numerical solution sampled at
cell centroids as a function of the radius $r=\|\mathbf x\|$, while the
red marker highlights the discrete boundary value at $r=R$ (the cut-cell
boundary/interfacial unknown used to impose Eq.~\eqref{eq:ball_robin}). The close
agreement with the analytical radial series Eq.~\eqref{eq:ball_series} confirms
that the moment-weighted boundary closure remains accurate in 3D and that the
boundary degree of freedom provides a consistent trace of the solution at
$\Gamma$.

\paragraph{Neumann limit and global conservation.}
To verify global conservation in the pure Neumann limit, we set
\(
\partial_n \phi = 0 \text{ on } \Gamma .
\)
The analytical configuration is chosen such that the domain-integrated
quantity is constant in time. On the discrete level, we monitor the conserved
integral and its relative drift,
\[
I^h(t) = \sum_{\text{active cells } i} V_i\,\phi_i(t), \quad
\delta I^h(t)=\frac{I^h(t)-I^h(0)}{I^h(0)}.
\]
Using the same time-integration parameters as in the Robin tests (midpoint
rule, $t_f=0.1$), $I^h(t)$ remains constant up to machine precision: the maximum
drift over $[0,t_f]$ is $|\delta I^h|_{\max}=3.1\times 10^{-17}$, and at the final
time we obtain $\delta I^h(t_f)=2.8\times 10^{-17}$.

\subsection{Two-phase validation}
\label{sec:diph_validation}

We now validate the two-phase formulation on several diffusion
problem.

\subsubsection{One-dimensional two-phase unsteady diffusion with homothetic jump}
\label{sec:diph_1d_jump}

We consider a one-dimensional two-phase parabolic diffusion problem as a check of the jump enforcement in a minimal setting. The domain is the
interval $\Omega=[0,L]$ split by a fixed interface at $x=x_{\mathrm{int}}$ into
\begin{equation}
\Omega^-=\{x\in[0,L]:\ x<x_{\mathrm{int}}\},
\qquad
\Omega^+=\{x\in[0,L]:\ x>x_{\mathrm{int}}\},
\label{eq:1d_split}
\end{equation}
where $x_{\mathrm{int}}$ is chosen so as not to coincide with a Cartesian face,
so that the interface cuts a control volume.

In each phase we solve the unsteady diffusion equation with equal diffusivities
$K^- = K^+ = K$,
\begin{equation}
\partial_t \phi^\pm = K\,\partial_{xx}\phi^\pm
\qquad \text{in }\Omega^\pm,
\label{eq:1d_parabolic}
\end{equation}
At the interface, we enforce continuity of diffusive flux together with a
homothetic jump in the field:
\begin{equation}
\llbracket K\,\partial_x \phi \rrbracket = 0,
\qquad
\phi^+ - \lambda\,\phi^- = 0
\qquad \text{at }x=x_{\mathrm{int}},
\label{eq:1d_interface}
\end{equation}
where $\lambda$ is prescribed.
A self-similar solution is obtained by introducing the similarity variables,
\[
\eta^\pm(x,t)=\frac{x-x_{\mathrm{int}}}{2\sqrt{K t}},
\]
and defining the constant
\begin{equation}
A
=
-\frac{\lambda}{1+\lambda}.
\label{eq:1d_prefactor_A}
\end{equation}
The exact
solution to the problem Eqs.~\eqref{eq:1d_parabolic}-\eqref{eq:1d_interface} for all
$t>0$ is then given by
\begin{align}
\phi^{-,ex}(x,t)
&=
A\Bigl[\operatorname{erfc}\!\bigl(\eta^-(x,t)\bigr)-2\Bigr],
\qquad x<x_{\mathrm{int}},
\label{eq:1d_phi_minus}\\[0.25em]
\phi^{+,ex}(x,t)
&=
A\,\operatorname{erfc}\!\bigl(\eta^+(x,t)\bigr)+1,
\qquad x>x_{\mathrm{int}}.
\label{eq:1d_phi_plus}
\end{align}
The outer boundary conditions are then prescribed as time-dependent Dirichlet
condition taken from Eqs.~\eqref{eq:1d_phi_minus}-\eqref{eq:1d_phi_plus}:
\[
\phi^-(0,t)=\phi^{-,\rm ex}(0,t),\qquad
\phi^+(L,t)=\phi^{+,\rm ex}(L,t),
\]
and the initial condition too. In the results below, we take $\lambda=100$ and
$K^-{=}K^+$.

Tab.~\ref{tab:1d_diph_lambda} reports $L^2$ errors at final time $t_f$
over regular, cut and all active cells. Pairwise convergence orders are also provided.

\begin{table}[h!]
\centering
\setlength{\tabcolsep}{6pt}
\renewcommand{\arraystretch}{1.2}
\begin{tabular}{c c | c c c | c c c}
\hline
$\lambda$ & $h$ &
$\|e(t_f)\|_{2,\mathrm{reg}}$ &
$\|e(t_f)\|_{2,\mathrm{cut}}$ &
$\|e(t_f)\|_{2,\mathrm{all}}$ &
$p_{\mathrm{reg}}$ & $p_{\mathrm{cut}}$ & $p_{\mathrm{all}}$ \\
\hline
100 & 2.0     & $2.16e{-1}$ & $2.90e{-1}$ & $3.61e{-1}$ & -   & -   & -   \\
100 & 1.0     & $1.85e{-1}$ & $1.29e{-1}$ & $2.25e{-1}$ & 0.22 & 1.17 & 0.68 \\
100 & 0.5     & $4.42e{-2}$ & $3.62e{-2}$ & $5.71e{-2}$ & 2.06 & 1.84 & 1.98 \\
100 & 0.25    & $2.74e{-2}$ & $1.14e{-3}$ & $2.74e{-2}$ & 0.69 & 4.99 & 1.06 \\
100 & 0.125   & $3.45e{-3}$ & $1.21e{-4}$ & $3.45e{-3}$ & 2.99 & 3.24 & 2.99 \\
100 & 0.0625  & $1.88e{-3}$ & $2.44e{-5}$ & $1.88e{-3}$ & 0.88 & 2.31 & 0.88 \\
100 & 0.03125 & $2.16e{-4}$ & $9.07e{-7}$ & $2.16e{-4}$ & 3.12 & 4.75 & 3.12 \\
\hline
fit & -     & -          & -          & -          & 2.00 & 3.53 & 2.00 \\
\hline
\end{tabular}
\caption{One-dimensional two-phase unsteady diffusion with homothetic jump
$\lambda=100$: $L^2$ errors at $t_f$ over regular, cut and all cells,
and corresponding convergence orders.}
\label{tab:1d_diph_lambda}
\end{table}

In addition to the stiff case $\lambda=100$ reported below, we also consider a
small sweep of homothetic jumps (e.g.\ $\lambda\in\{0.1,1,10,100\}$) in order to
visualize the interface enforcement. As $\lambda$ varies, the analytical
solution Eqs.~\eqref{eq:1d_phi_minus}-\eqref{eq:1d_phi_plus} transitions from a
weak to a strong discontinuity in $\phi$ across $x_{\mathrm{int}}$, while the
flux continuity constraint in Eq.~\eqref{eq:1d_interface} remains unchanged. This
makes the test particularly sensitive to the discrete coupling: the method must
simultaneously reproduce the prescribed jump amplitude and maintain consistent
normal gradients on both sides of an interface that cuts a control volume.

\begin{figure}[!h]
  \centering
  \includegraphics[width=0.86\linewidth]{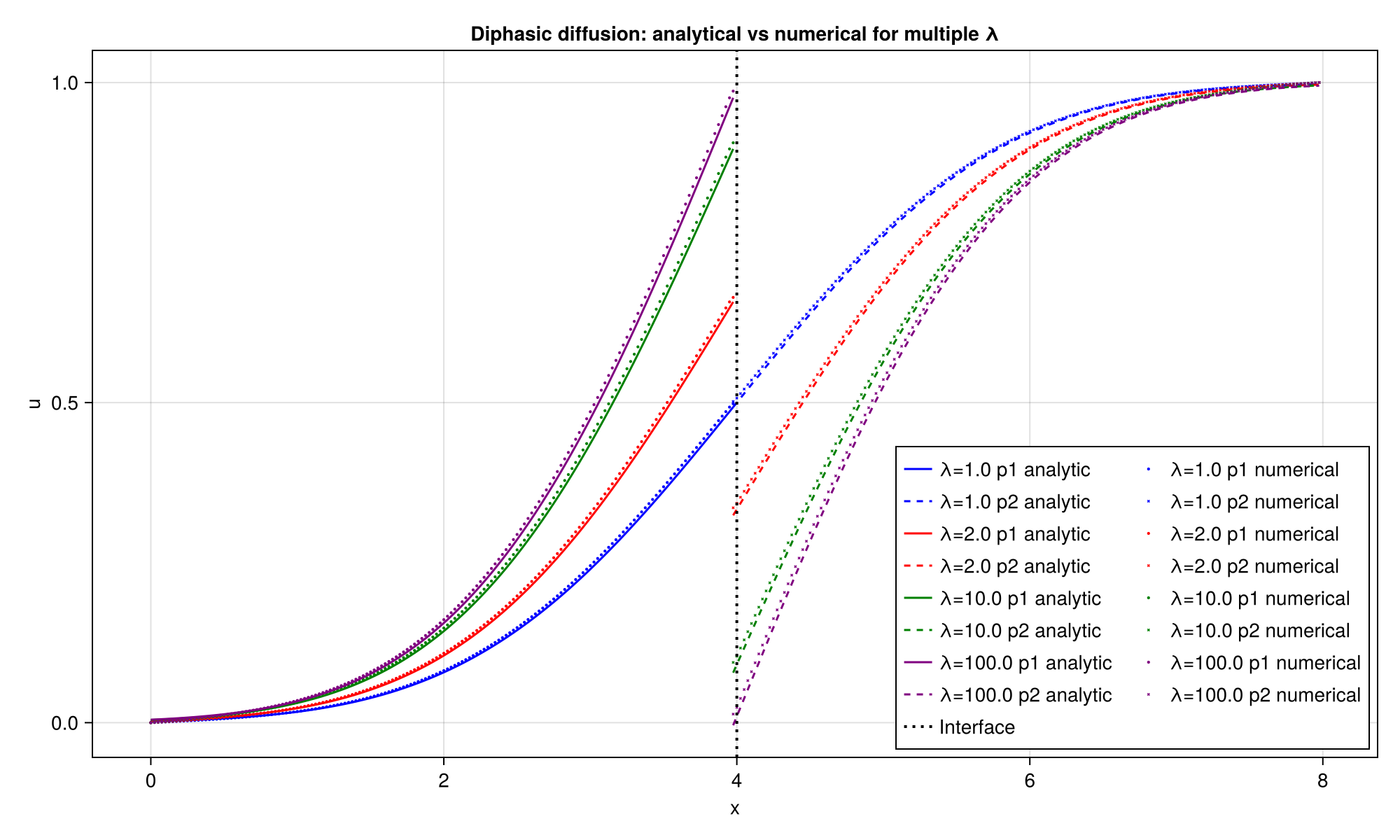}
  \caption{One-dimensional two-phase unsteady diffusion with homothetic jump:
  numerical solution (markers) and analytical profiles (lines) at the
  final time $t_f$ for several jump ratios $\lambda$ (e.g.\ $\lambda=0.1,1,10,100$).
  The vertical dashed line indicates the interface position $x_{\mathrm{int}}$
  (chosen inside a cut cell : $x_{int}=4.0+1e^{-2}$). The correct jump amplitude $\phi^+(x_{\mathrm{int}},t_f)
  = \lambda\,\phi^-(x_{\mathrm{int}},t_f)$ is recovered while preserving flux
  continuity across the interface.}
  \label{fig:1d_jump_lambda_sweep}
\end{figure}
For all $\lambda$ values tested, the numerical profiles match the analytical
solution near the interface and the discrepancy remains localized to a small
neighborhood of the cut cell, confirming that the coupling is both sharp and
stable even for large jump ratios.

\subsubsection{Two-phase unsteady diffusion with a circular interface}
\label{sec:diph_circle_unsteady}
The computational domain is the square cavity
\(
\Omega = [0,L]^2, L=8,
\)
containing a circular interface $\Gamma$ of radius $R_0=2$ centered at
$\mathbf x_c=(4,4)$. The phases are defined by
\[
\Omega^+ = \{\mathbf x\in\Omega:\ r<R_0\}, \quad
\Omega^- = \{\mathbf x\in\Omega:\ r>R_0\}, \quad
r=\|\mathbf x-\mathbf x_c\|.
\]
In each phase, we solve the
unsteady diffusion equation with homogeneous Neumann conditions on the outer
boundary $\partial\Omega$. At the interface, we consider the continuity case
$\lambda=1.0$, i.e.\ both the field and the normal diffusive flux are continuous
across $\Gamma$. The initial data are chosen as $\phi^+(\mathbf x,0)=\phi_0$ and
$\phi^-(\mathbf x,0)=0$.

For this radially symmetric configuration, an analytical reference solution is
available in polar coordinates $(r,\theta)$ in the form of integral
representations involving Bessel functions.
Denoting $K=\sqrt{K^+/K^-}$ and letting $\Phi(u)$ and $\Psi(u)$ be the auxiliary
functions defined below, the solution reads
\begin{align}
\phi^+(r,t)
&= \frac{4\,\phi_0\,K^+\,(K^-)^2\,}{\pi^2\,R_0}
    \int_{0}^{\infty}
      e^{-K^+\,u^2\,t}
      \,\frac{J_0(u\,r)\,J_1(u\,R_0)}
             {u^2\bigl[\Phi(u)^2+\Psi(u)^2\bigr]}
    \,\mathrm{d}u,
\label{eq:phi_plus_analytical}\\
\phi^-(r,t)
&= \frac{2\,\phi_0\,K^+\,\sqrt{K^-}\,}{\pi}
    \int_{0}^{\infty}
      e^{-K^+\,u^2\,t}
      \,\frac{J_1(u\,R_0)\,\bigl[J_0\bigl(K\,u\,r\bigr)\,\Phi(u)
                             -Y_0\bigl(K\,u\,r\bigr)\,\Psi(u)\bigr]}
             {u\bigl[\Phi(u)^2+\Psi(u)^2\bigr]}
    \,\mathrm{d}u,
\label{eq:phi_minus_analytical}
\end{align}
and
\begin{align*}
  \Phi(u)
  &= K^+\,\sqrt{K^-}\;J_1(R_0u)\;Y_0\bigl(K\,R_0u\bigr)
    \;-K^-\,\sqrt{D^+}\;J_0(R_0u)\;Y_1\bigl(K\,R_0u\bigr),\\
  \Psi(u)
  &= K^+\,\sqrt{K^-}\;J_1(R_0u)\;J_0\bigl(K\,R_0u\bigr)
    \;-K^-\,\sqrt{D^+}\;J_0(R_0u)\;J_1\bigl(K\,R_0u\bigr).
\end{align*}

Tab.~\ref{tab:2ph-2D-circ-final} reports the $L^2$ errors at the final
time $t_f$ over regular cells, cut cells and all active cells in both phases.
Several features are worth noting.

\begin{table}[!h]
\centering
\setlength{\tabcolsep}{6pt}
\renewcommand{\arraystretch}{1.2}
\begin{tabular}{c c | c c c | c c c}
\hline
$h$ &
$N_{diam}$ &
$\|e(t_f)\|_{2,\mathrm{reg}}$ &
$\|e(t_f)\|_{2,\mathrm{cut}}$ &
$\|e(t_f)\|_{2,\mathrm{all}}$ &
$p_{\mathrm{reg}}$ &
$p_{\mathrm{cut}}$ &
$p_{\mathrm{all}}$
\\
\hline
2.0    & $ 1$  & $4.253e{-1}$ & $5.620e{-1}$ & $7.048e{-1}$ &  -   &  -   &  -   \\
1.0    & $ 2$  & $2.363e{-1}$ & $2.536e{-1}$ & $3.467e{-1}$ & 0.85  & 1.15  & 1.02  \\
0.5    & $ 8$  & $8.458e{-2}$ & $4.781e{-2}$ & $9.716e{-2}$ & 1.48  & 2.41  & 1.84  \\
0.25   & $12$  & $4.358e{-2}$ & $7.874e{-3}$ & $4.429e{-2}$ & 0.96  & 2.60  & 1.13  \\
0.125  & $24$  & $7.989e{-3}$ & $1.190e{-3}$ & $8.077e{-3}$ & 2.45  & 2.73  & 2.45  \\
0.0625 & $48$  & $3.215e{-3}$ & $2.951e{-4}$ & $3.229e{-3}$ & 1.31  & 2.01  & 1.32  \\
\hline
fit    & -    & - & - & - & 1.88 & 2.37 & 1.89 \\
\hline
\end{tabular}
\caption{$L^2$ error at final time $t_f$ for the 2D two-phase diffusion
test with a circular interface.}
\label{tab:2ph-2D-circ-final}
\end{table}

For the coarsest resolutions ($h=2$ and $h=1$), the circular interface is
severely under-resolved and only a handful of cut cells represent the coupling
between phases. Even though, the solver remains robust with at most 1 cell per domain diameter.
Once the interface is resolved by $O(10)$-$O(50)$ cells across the relevant
length scale (here, from $h=0.5$ down to $h=0.0625$), the global error approaches
a nearly second-order trend: the fitted rates are
$p_{\mathrm{all}}\approx 1.89$ and $p_{\mathrm{reg}}\approx 1.88$. The global
order is slightly below $2$ because the solution remains influenced by the
embedded interface region, where the local stencil is modified and the regular
Cartesian cancellation is partially lost. Importantly, the convergence of the
all-cells metric indicates that the two-fluid coupling does not spoil the
accuracy of the bulk discretization. The cut-cell error decreases faster than the bulk error in this test, with a fitted rate $p_{\mathrm{cut}}\approx 2.37$. This comparatively high rate is
consistent with the fact that the coupling is enforced strongly through
interfacial unknowns and moment-weighted operators, so that the dominant
interface error can drop rapidly once smooth geometric under-resolution is removed. Some pairwise rates (e.g.\ the drop in $p_{\mathrm{all}}$ between $h=0.5$ and
$h=0.25$) reflect the discrete change in the cut-cell topology as the interface
position change relative to the background grid. As $h$ changes, the set of cut cells and
their aperture patterns are altered, which can temporarily bias the balance
between bulk and interface contributions in the error norm. The convergence
trend is nonetheless clear on the finest grids.

The temporal evolution of the interfacial diffusive flux provides an additional
consistency check for the two-phase coupling. At each time step we evaluate
\[
Q^+(t)=\int_{\Gamma} K^+\,\nabla\phi^+\cdot\mathbf n\,\mathrm dS,
\]
which measures the instantaneous transfer rate across the interface.

\begin{figure}[!h]
    \centering
    \includegraphics[width=0.75\linewidth]{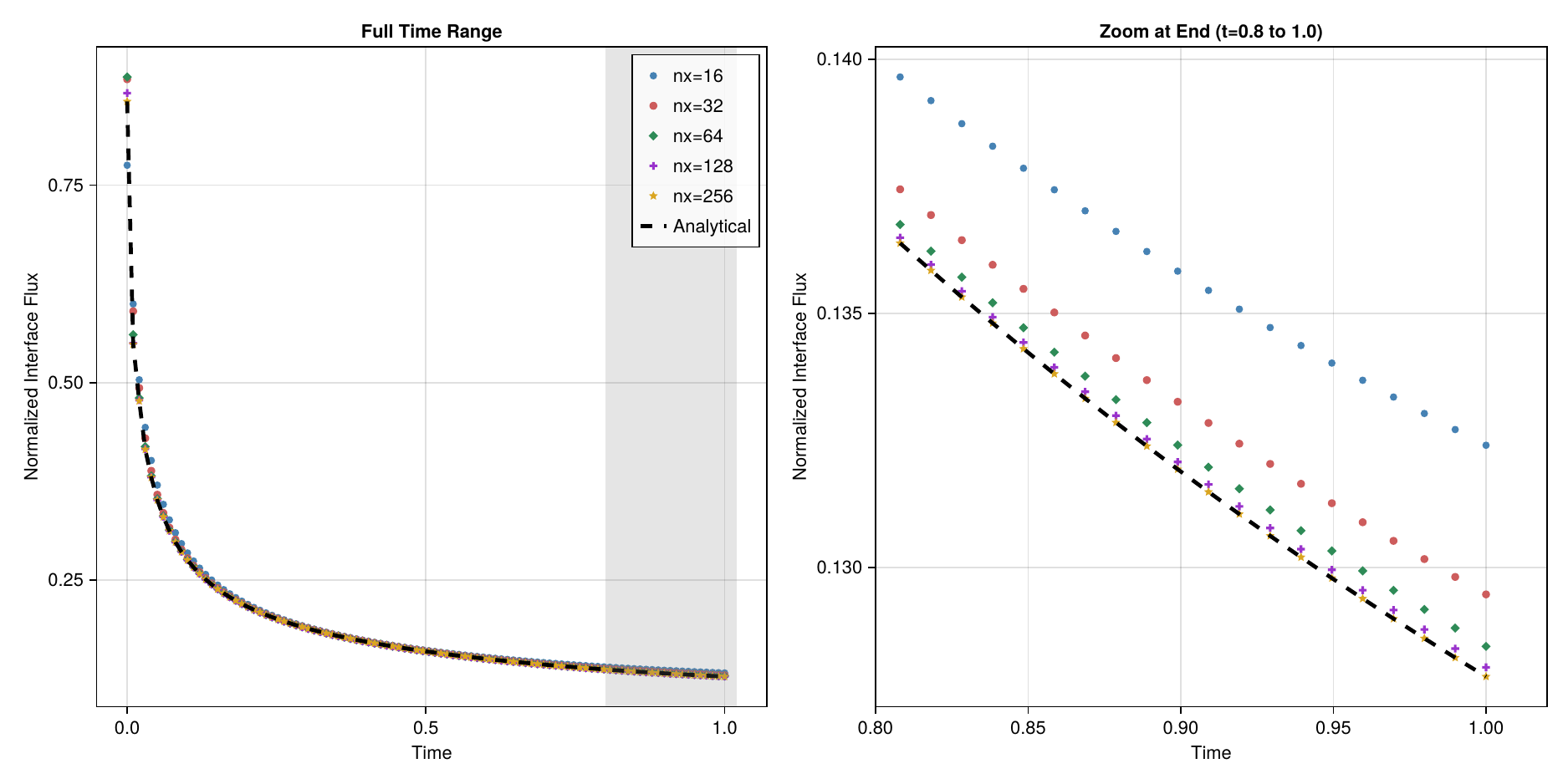}
    \caption{Two-phase unsteady diffusion with a circular interface: Normalized interfacial flux across the interface}
    \label{fig:diph_circle}
\end{figure}

The temporal evolution of the interfacial diffusive flux
provides a stringent check of the two-phase coupling. Tab.~\ref{tab:bi_flux_rel}
shows that the final flux value $Q^+(t_f)$ converges monotonically toward
its reference value as the mesh is refined, with the relative error decreasing
from $3.56\times 10^{-2}$ at $h=0.5$ to $1.42\times 10^{-3}$ at $h=0.0625$.
The corresponding rates are close second order over this range,
indicating that both the enforcement of the interfacial condition and the reconstruction of normal gradients in
cut cells remain accurate. Fig.~\ref{fig:diph_circle} further illustrates that
the normalized interfacial flux history collapses under refinement:
transient discrepancies on coarse grids diminish confirming that the method captures the
time-dependent exchange across the embedded interface.
\begin{table}[h!]
\centering
\renewcommand{\arraystretch}{1.2}
\begin{tabular}{c|c|c}
\hline
Resolution $h$ &$Q^+(t_f)$ & Relative error \\
\hline
0.5    & 0.1324020839 & $3.560e{-2}$ \\
0.25   & 0.1294676902 & $1.265e{-2}$ \\
0.125  & 0.1284444328 & $4.642e{-3}$ \\
0.0625 & 0.1280328469 & $1.423e{-3}$ \\
\hline
\end{tabular}
\caption{Final interfacial flux diagnostic $Q^+(t_f)$ and relative error for the
two-phase circular test.}
\label{tab:bi_flux_rel}
\end{table}

\subsubsection{Two-phase unsteady diffusion across a sphere}
\label{sec:val_diph_3D_sphere_brown}

We conclude the validation with a three-dimensional two-phase test in
a spherical geometry, for which a closed analytical reference is available in
integral form \citep{brown_diffusion_1965}. This benchmark simultaneously
checks  the two-fluid coupling across a curved interface through scalar
continuity and flux continuity and  the robustness of the 3D
moment integration 
required by the cut-cell discretization.

Let $\mathbf x_c=(2,2,2)$ and $R_0=1$. The interface is the sphere
\(
\Gamma=\{\mathbf x:\ r=R_0\},  r=\|\mathbf x-\mathbf x_c\|,
\)
splitting $\Omega_{\rm box}=[0,4]^3$ into the inner phase
$\Omega^-=\{r<R_0\}$ and the outer phase $\Omega^+=\{r>R_0\}$. 
We consider unsteady diffusion in the two phases $\Omega^\pm$, with constant
diffusivities $D^-$ and $D^+$. Across the interface $\Gamma$ at
$r=R_0$, we enforce continuity of the scalar field together with continuity of
the normal diffusive flux.

The analytical solution corresponds to the classical configuration of an
initially uniform inner concentration and vanishing outer concentration,
\begin{equation}
\phi^-(r,0)=\phi_0, \quad 0\le r<R_0,
\qquad
\phi^+(r,0)=0, \quad r>R_0,
\label{eq:brown_ic}
\end{equation}
together with decay at infinity, $\phi^+(r,t)\to 0$ as $r\to\infty$.
In our bounded box, Dirichlet condition on $\partial\Omega_{\rm box}$ are taken from
the analytical solution evaluated at the boundary, so that the truncation does
not pollute the error assessment.

The solution is radially symmetric, $\phi^\pm(\mathbf x,t)=\Phi^\pm(r,t)$.
Following \cite{brown_diffusion_1965}, we introduce the auxiliary parameters
\begin{equation}
\sigma = \sqrt{D^+/D^-},
\qquad
Q = \frac{D^+}{D^-}\,\sigma,
\qquad
L = \frac{D^+-D^-}{D^-},
\label{eq:brown_params}
\end{equation}
and the denominator
\begin{equation}
\mathcal D(u)
=
\bigl(u\cos u + L \sin u\bigr)^2
+
\bigl(Q\,u\sin u\bigr)^2.
\label{eq:brown_denom}
\end{equation}
For the inner phase ($0<r<R_0$), the exact solution reads
\begin{equation}
\Phi^-(r,t)
=
\frac{2Q\alpha \phi_0}{\pi r}
\int_0^{\infty}
\frac{\bigl(\sin u-u\cos u\bigr)\,
\sin\!\left(ur/R_0\right)}
{\mathcal D(u)}\;
\exp\!\left(-D^- u^2 t/R_0^2\right)\,{\rm d}u.
\label{eq:brown_core}
\end{equation}
For the outer phase ($r>R_0$), introduce
\begin{equation}
\mathcal F(u,r)
=
\bigl(u\cos u + L \sin u\bigr)\,
\sin\!\left(u(r-R_0)/\sigma R_0\right)
+
Q\,u\sin u\,
\cos\!\left(u(r-R_0)/\sigma R_0\right),
\label{eq:brown_F}
\end{equation}
and the exact solution is
\begin{equation}
\Phi^+(r,t)
=
\frac{2\alpha \phi_0}{\pi r}
\int_0^{\infty}
\frac{\bigl(\sin u-u\cos u\bigr)\,
\mathcal F(u,r)}
{u\,\mathcal D(u)}\;
\exp\!\left(-D^- u^2 t/R_0^2\right)\,{\rm d}u.
\label{eq:brown_shell}
\end{equation}
In practice, the integrals Eqs.~\eqref{eq:brown_core}-\eqref{eq:brown_shell} are
evaluated by adaptive quadrature and truncated at a large cutoff. Time integration uses a midpoint time scheme, with a time step $\Delta t = 0.25 \min(\Delta x, \Delta y, \Delta z)$.

Tab.~\ref{tab:diph_3D_sphere_brown} reports the $L^2$ errors at the final time
$t_f$ over all active cells, regular cells and cut cells. The systematic decay
of these three measures under mesh refinement confirms that the two-fluid
two-phase coupling remains consistent in three dimensions. This benchmark is
also a stringent stress test for the geometric layer: a spherical interface
generates a broad spectrum of cut-cell patterns, including very small trimmed
control volumes, strongly skewed face apertures and highly curved interface
patches. The convergence observed in Tab.~\ref{tab:diph_3D_sphere_brown}
therefore indicates that the three-dimensional moment integration (phase
volumes, face apertures and interface measures/centroids) is sufficiently
accurate and robust to sustain stable discrete flux balances and sharp
enforcement of the interface conditions. The pairwise orders are not uniform across all refinement steps. For coarse
meshes, the sphere is under-resolved and changes in cut-cell topology can
temporarily shift the relative weight of bulk and interface contributions,
leading to pre-asymptotic behavior. As the resolution increases, the regular-
cell and global errors display a clear rapid decay, whereas the cut-cell error
may exhibit apparent superconvergence.

\begin{table}[!h]
\centering
\setlength{\tabcolsep}{6pt}
\renewcommand{\arraystretch}{1.2}
\begin{tabular}{c c | c c c | c c c}
\hline
$h$ & $N$ &
$\|e(t_f)\|_{2,\mathrm{all}}^{\rm rel}$ &
$\|e(t_f)\|_{2,\mathrm{reg}}^{\rm rel}$ &
$\|e(t_f)\|_{2,\mathrm{cut}}^{\rm rel}$ &
$p_{\mathrm{all}}$ & $p_{\mathrm{reg}}$ & $p_{\mathrm{cut}}$ \\
\hline
0.5         & 4  & $4.845e{-2}$ & $2.588e{-2}$ & $4.096e{-2}$ &  -   &  -   &  -   \\
0.25        & 8  & $9.087e{-3}$ & $8.116e{-3}$ & $4.085e{-3}$ & 2.41 & 1.67 & 3.33 \\
0.12903226  & 16 & $2.748e{-3}$ & $2.441e{-3}$ & $1.261e{-3}$ & 1.81 & 1.82 & 1.78 \\
0.08333333  & 24 & $6.286e{-4}$ & $6.251e{-4}$ & $6.573e{-5}$ & 3.37 & 3.12 & 6.76 \\
0.04166667  & 48 & $1.771e{-4}$ & $1.676e{-4}$ & $2.167e{-5}$ & 1.83 & 1.90 & 1.60 \\
\hline
fit  & - & - & - & - & 2.28 & 2.07 & 3.13 \\
\hline
\end{tabular}
\caption{Two-phase unsteady diffusion in a sphere : $L^2$
errors at final time $t_f$ over all, regular and cut cells, with pairwise
convergence orders.}
\label{tab:diph_3D_sphere_brown}
\end{table}

\section{Conclusions and outlook}
\label{sec:conclusion}

This work introduced a conservative Cartesian cut-cell finite-volume method
for sharp-interface diffusion problems in static embedded geometries,
in both single-phase and two-phase configurations. The formulation relies on a
reduced geometric description: phase-restricted volumes, face apertures and a
small set of low-order geometric moments that are sufficient to assemble compact
stencils for discrete divergence and gradient operators. Boundary conditions (Dirichlet, Neumann and Robin) and interface jumps (normal diffusive flux and value) are enforced strongly in the algebraic system. A distinctive feature of the two-fluid approach is that interfacial unknowns are explicit degrees of freedom, providing direct access to interface values
without post-processing or interpolation.

The numerical validation in static domains demonstrates accuracy, conservation,
and robustness across a range of embedded geometries. Single-phase
elliptic tests on non-convex domains recover near second-order
convergence in $L^2$ norm, with cut-cell errors converging at comparable rates. The discrete solution remains within the expected physical range, with no overshoot or undershoot observed as the grid is refined. For Robin boundary condition, the method maintains robust convergence
in both $L^2$ and $H^1$ norms and preserves a faithful reconstruction of boundary
unknowns. In three dimensions, the method extends without
modification and remains effective even at very low resolutions, when the
embedded body is represented by only a few cells. In the Neumann limit, the
discrete mass is preserved to machine precision, highlighting exact flux
conservation at the discrete level. Two-phase validations confirm that sharp coupling across interfaces does not compromise robustness or accuracy. The method exhibits near second-order convergence in global norms and appropriate behavior in cut cells, while the interfacial flux converges rapidly under refinement. Even with a strong homothetic jump, the coupling remains robust. Both the jump enforcement and the flux continuity remain well behaved as the interface cuts control volumes at arbitrary locations. Overall, the results show that the proposed method provides a reliable and geometrically flexible building block for interface-resolved diffusion with compact conservative operators constructed from a reduced set of moments, direct discrete enforcement of boundary and interface conditions, robust behavior in severely under-resolved regimes and consistent accuracy in both mono- and two-phase configurations.

While this article focuses on fixed interfaces, the formulation is designed as a
foundation for more general sharp-interface models. Natural extensions include
prescribed-motion geometries, free-boundary diffusion
problems in which the interface motion is determined by interfacial balances,
and Stefan-type phase-change models coupling diffusion to latent-heat
conditions. Beyond pure diffusion, the same cut-cell machinery can be combined
with sharp-interface incompressible flow solvers to address coupled
advection-diffusion and two-phase transport. These developments will leverage
the core properties established here ; local conservation, accurate interface
coupling and robustness on complex embedded geometries, to tackle fully coupled
multiphysics problems with dynamically evolving interfaces.

\newpage
\appendix

\include{appendix/table_static}

\include{appendix/dem_inver}
\include{appendix/1ddev}

\bibliographystyle{plain}
\bibliography{references}

\end{document}

%% file: appendix/table_static.tex
\section{Summary table for Static quantities notations}
\begin{table}[h!]
\centering
\renewcommand{\arraystretch}{1.3}
\begin{tabular}{lll}
\hline
\textbf{Quantity} & \textbf{Definition using $\langle \cdot , \cdot \rangle$} & \textbf{Equivalent integral form} \\[3pt]
\hline
Cell volume &
$V_{i,j}^{\pm}(t) = \langle \Omega_{i,j}^{\pm}(t), 1 \rangle$ &
$\displaystyle \int_{\Omega_{i,j}^{\pm}(t)} 1\,\mathrm{d}V$ \\[4pt]

Cell centroid &
$x_{i,j}^{\pm}(t) = \dfrac{\langle \Omega_{i,j}^{\pm}(t), x \rangle}{\langle \Omega_{i,j}^{\pm}(t), 1 \rangle}$,\quad
\text{eq for y}&
$\displaystyle \frac{\int_{\Omega_{i,j}^{\pm}(t)} x\,\mathrm{d}V}{\int_{\Omega_{i,j}^{\pm}(t)} \mathrm{d}V}$,
\text{eq for y} \\[6pt]

Face area (x–normal) &
$A_{i-\sfrac{1}{2},j}^{1\pm}(t) = \langle \Sigma_j^{1\pm}(t, x_{i-\sfrac{1}{2}}), 1 \rangle$ &
$\displaystyle \int_{\Sigma_j^{1\pm}(t, x_{i-\sfrac{1}{2}})} 1\,\mathrm{d}S$ \\[4pt]

Face area (y–normal) &
$A_{i,j-\sfrac{1}{2}}^{2\pm}(t) = \langle \Sigma_i^{2\pm}(t, y_{j-\sfrac{1}{2}}), 1 \rangle$ &
$\displaystyle \int_{\Sigma_i^{2\pm}(t, y_{j-\sfrac{1}{2}})} 1\,\mathrm{d}S$ \\[4pt]

Centroidal face area (x–normal) &
$B_{i,j}^{1\pm}(t) = \langle \Sigma_j^{1\pm}(t, x_{i,j}^{\pm}(t)), 1 \rangle$ &
$\displaystyle \int_{\Sigma_j^{1\pm}(t, x_{i,j}^{\pm}(t))} 1\,\mathrm{d}S$ \\[4pt]

Centroidal face area (y–normal) &
$B_{i,j}^{2\pm}(t) = \langle \Sigma_i^{2\pm}(t, y_{i,j}^{\pm}(t)), 1 \rangle$ &
$\displaystyle \int_{\Sigma_i^{2\pm}(t, y_{i,j}^{\pm}(t))} 1\,\mathrm{d}S$ \\[4pt]

Staggered volume (x–direction) &
$W_{i-\sfrac{1}{2},j}^{1\pm}(t) = \langle \Omega_j^{1\pm}(t, \left] x_{i-1,j}^{\pm}, x_{i,j}^{\pm} \right[), 1 \rangle$ &
$\displaystyle \int_{\Omega_j^{1\pm}(t; x_{i-1,j}^{\pm}, x_{i,j}^{\pm})} 1\,\mathrm{d}V$ \\[4pt]

Staggered volume (y–direction) &
$W_{i,j-\sfrac{1}{2}}^{2\pm}(t) = \langle \Omega_i^{2\pm}(t, \left] y_{i,j-1}^{\pm}, y_{i,j}^{\pm} \right[), 1 \rangle$ &
$\displaystyle \int_{\Omega_i^{2\pm}(t; y_{i,j-1}^{\pm}, y_{i,j}^{\pm})} 1\,\mathrm{d}V$ \\[3pt]
\hline\hline

Bulk variable average &
$\displaystyle \Phi_{i,j}^{\omega\pm}(t) = 
\frac{\langle \Omega_{i,j}^{\pm}(t), \phi^{\pm}(t) \rangle}{V_{i,j}^{\pm}(t)}$ &
$\displaystyle \frac{\int_{\Omega_{i,j}^{\pm}(t)} \phi^{\pm}(t)\,\mathrm{d}V}
{\int_{\Omega_{i,j}^{\pm}(t)} \mathrm{d}V}$ \\[4pt]

Interface variable average &
$\displaystyle \Phi_{i,j}^{\gamma\pm}(t) =
\frac{\langle \Gamma_{i,j}(t), \phi^{\pm}(t) \rangle}{\langle \Gamma_{i,j}(t), 1 \rangle}$ &
$\displaystyle \frac{\int_{\Gamma_{i,j}(t)} \phi^{\pm}(t)\,\mathrm{d}S}
{\int_{\Gamma_{i,j}(t)} \mathrm{d}S}$ \\[4pt]

Flux component (x–direction) &
$\displaystyle Q_{i-\sfrac{1}{2},j}^{1\pm}(t) =
\frac{\langle \Sigma_j^{1\pm}(x_{i-\sfrac{1}{2}}), q^{1\pm}(t) \rangle}
{A_{i-\sfrac{1}{2},j}^{1\pm}}$ &
$\displaystyle \frac{\int_{\Sigma_j^{1\pm}(x_{i-\sfrac{1}{2}})} q^{1\pm}(t)\,\mathrm{d}S}
{\int_{\Sigma_j^{1\pm}(x_{i-\sfrac{1}{2}})} \mathrm{d}S}$ \\[4pt]

Source term average &
$\displaystyle R_{i,j}^{\pm}(t) =
\frac{\langle \Omega_{i,j}^{\pm}(t), r_{i,j}^{\pm}(t) \rangle}
{V_{i,j}^{\pm}(t)}$ &
$\displaystyle \frac{\int_{\Omega_{i,j}^{\pm}(t)} r_{i,j}^{\pm}(t)\,\mathrm{d}V}
{\int_{\Omega_{i,j}^{\pm}(t)} \mathrm{d}V}$ \\[4pt]

Interface field average &
$\displaystyle F_{i,j}(t) =
\frac{\langle \Gamma_{i,j}(t), f(t) \rangle}{\langle \Gamma_{i,j}(t), 1 \rangle}$ &
$\displaystyle \frac{\int_{\Gamma_{i,j}(t)} f(t)\,\mathrm{d}S}
{\int_{\Gamma_{i,j}(t)} \mathrm{d}S}$ \\[3pt]
\hline
\end{tabular}
\caption{Summary of geometric and semi–discrete quantities expressed using the generic integral notation $\langle \Xi, f \rangle$.}
\label{tab:geom_and_sd}
\end{table}

%% file: appendix/dem_inver.tex
\section{Invertibility of the interfacial augmented system}

We rewrite the $4\times4$ block system from \ref{sec:discrete} as
\[
\underbrace{\begin{bmatrix}
    M^- + \theta L^{\omega\omega-} & 0 & L^{\omega\gamma-} & 0 \\
    0 & M^+ + \theta L^{\omega\omega+} & 0 & L^{\omega\gamma+} \\
    -\theta L^{\gamma\omega-} & \theta L^{\gamma\omega+} & -L^{\gamma\gamma-} & L^{\gamma\gamma+} \\
    0 & 0 & -\lambda I & I
\end{bmatrix}}_{\displaystyle \mathcal M}
\begin{bmatrix}
\Phi^{\omega-}_{n+1}\\[2pt]
\Phi^{\omega+}_{n+1}\\[2pt]
\Phi^{\gamma-}_{n+\theta}\\[2pt]
\Phi^{\gamma+}_{n+\theta}
\end{bmatrix}
=\ \text{RHS}.
\]
For compactness, set
\[
\begin{aligned}
&A^-:=M^-+\theta L^{\omega\omega-},\qquad &&A^+:=M^++\theta L^{\omega\omega+},\\
&C^-:=L^{\omega\gamma-},\qquad &&C^+:=L^{\omega\gamma+},\\
&B^-:=L^{\gamma\omega-},\qquad &&B^+:=L^{\gamma\omega+},\\
&D^-:=L^{\gamma\gamma-},\qquad &&D^+:=L^{\gamma\gamma+}.
\end{aligned}
\]
Group bulk and interface unknowns as
\[
\Phi^\omega:=\begin{bmatrix}\Phi^{\omega-}\\ \Phi^{\omega+}\end{bmatrix},\qquad
\Phi^\gamma:=\begin{bmatrix}\Phi^{\gamma-}\\ \Phi^{\gamma+}\end{bmatrix}.
\]
Then
\[
\mathcal M=\begin{bmatrix}A & C\\[2pt] B & \widehat D\end{bmatrix},\qquad
A=\mathrm{diag}(A^-,A^+),\quad
C=\begin{bmatrix}C^-&0\\[2pt]0&C^+\end{bmatrix},\quad
B=\begin{bmatrix} - \theta B^- & \theta B^+\\[2pt] 0 & 0\end{bmatrix},\quad
\widehat D=\begin{bmatrix} -D^- & D^+\\[2pt] -\lambda I & I\end{bmatrix}.
\]

Standard assumptions leads to : 
\begin{enumerate}
    \item $M^\pm=\mathrm{diag}(V^\pm/\Delta t)\succ0$
    \item (ii) $L^{\omega\omega\pm}$ and $L^{\gamma\gamma\pm}$ are symmetric positive definite (SPD) whenever the face diffusion tensors are SPD (jacobians of symetric divergence and gradient operators)
\end{enumerate}
Hence $A^\pm$ and $A$ are SPD for any $\theta\in[0,1]$. In the typical discrete adjoint setting, $B^\pm=(C^\pm)^\top$ and $D^\pm=(D^\pm)^\top$.

Eliminate $\Phi^\gamma$ via the Schur complement of the interface block:
\[
\boxed{\,S_\omega \ :=\ A - C\,\widehat D^{-1}\,B\,}.
\]
The full system is invertible iff $\widehat D$ and $S_\omega$ are both invertible, because
$\det\mathcal M=\det\widehat D\cdot \det S_\omega$.

From the fourth row of $\widehat D$ we have $\Phi^{\gamma+}=\lambda\,\Phi^{\gamma-}$. Substituting into the third row yields
\[
\bigl(\lambda D^+ - D^-\bigr)\,\Phi^{\gamma-}
\;+\;\theta\,B^+\,\Phi^{\omega-}\;-\;\theta\,B^-\,\Phi^{\omega+}=0.
\]
Define the \emph{effective interface operator}
\[
\boxed{\,D_{\mathrm{eff}}:=\lambda D^+ - D^-\,}.
\]
Then
\[
\Phi^{\gamma-}
= -\,\theta\,D_{\mathrm{eff}}^{-1}\,\bigl(B^+\,\Phi^{\omega-}-B^-\,\Phi^{\omega+}\bigr),
\qquad
\Phi^{\gamma+}=\lambda\,\Phi^{\gamma-}.
\]
Consequently, $\widehat D$ is invertible iff $D_{\mathrm{eff}}$ is invertible. Since $D^\pm$ are SPD,
\[
D_{\mathrm{eff}}=\lambda D^+ - D^-
=(D^+)^{1/2}\Bigl[\lambda I - (D^+)^{-1/2}D^-(D^+)^{-1/2}\Bigr](D^+)^{1/2}.
\]
Thus $D_{\mathrm{eff}}\succ0$ whenever
\[
\boxed{\,\lambda > \lambda^{\max}_\Gamma,\qquad
\lambda^{\max}_\Gamma:=\lambda_{\max}\!\bigl((D^+)^{-1/2} D^- (D^+)^{-1/2}\bigr)\,.\,}
\]
Equivalently, for all $y\neq0$,
$y^\top D_{\mathrm{eff}} y
=\lambda\,y^\top D^+ y - y^\top D^- y
\ge (\lambda-\lambda^{\max}_\Gamma)\,y^\top D^+ y$.

For invertibility of $S_\omega$. Since $A^\pm\succ0$, there exists $\alpha>0$ such that
\[
\langle x,Ax\rangle \ \ge\ \alpha\,\|x\|^2\qquad\forall\,x=\begin{bmatrix}x^-\\ x^+\end{bmatrix}.
\]
Moreover,
\[
\bigl\|\,C\,D_{\mathrm{eff}}^{-1}B\,\bigr\|
\ \le\ \|C\|\,\|D_{\mathrm{eff}}^{-1}\|\,\|B\|.
\]
Hence, if
\[
\boxed{\,\theta\,\|C\|\,\|D_{\mathrm{eff}}^{-1}\|\,\|B\|\ <\ \alpha,}
\]
then for all $x\neq0$,
\[
\langle x,S_\omega x\rangle
=\langle x,Ax\rangle - \theta\,\langle x, C \widehat{D}^{-1} B x\rangle
\ \ge\ (\alpha-\theta\|C\|\,\|D_{\mathrm{eff}}^{-1}\|\,\|B\|)\,\|x\|^2\ >\ 0,
\]
Hence, $S_\omega$ is coercive and therefore invertible, since coercivity implies a trivial nullspace and a strictly positive Rayleigh quotient, i.e. all eigenvalues are positive and $S_\omega$ is SPD.

\emph{Remarks:} 
\begin{itemize}
    \item For $\theta=0$ (forward Euler), $S_\omega=A$ and is immediately invertible
    \item For $\theta>0$, the mass terms in $A^\pm$ increase $\alpha$, enlarging the parameter window.
\end{itemize}

With $M^\pm\succ0$ and $L^{\omega\omega\pm}$, $L^{\gamma\gamma\pm}$ are SPD then $D_{\mathrm{eff}}:=\lambda D^+- D^-$ is SPD and the bulk--only Schur complement
\[
S_\omega
=
\mathrm{diag}(A^-,A^+)\;-\;\theta
\begin{bmatrix} C^-\\ \lambda C^+\end{bmatrix}
D_{\mathrm{eff}}^{-1}
\begin{bmatrix} -B^- & \,B^+\end{bmatrix}
\]
is invertible. Consequently, the full matrix $\mathcal M$ is invertible.

%% file: appendix/1ddev.tex
\section{One‐dimensional two‐phase diffusion problem}
\label{appendix:1D_diffusion}

We consider three Cartesian cells of width \(\Delta x\) covering \(x\in[0,3\Delta x]\), with the middle cell cut by a sharp interface at \(x_\Gamma\in(\Delta x,2\Delta x)\).  Denote
\(
  \ell^- = x_\Gamma - \Delta x,\qquad
  \ell^+ = 2\Delta x - x_\Gamma,\qquad
  \ell^- + \ell^+ = \Delta x.
\)
At time level \(n\) we have four bulk averages
\(
  \Phi_1^n,\quad \Phi_{2}^{-,n},\quad \Phi_{2}^{+,n},\quad \Phi_3^n,
\)
and two interface‐side values
\(
  \Phi_{\Gamma}^{-,n},\quad \Phi_{\Gamma}^{+,n}.
\)

\paragraph{Finite‐volume balances}  
With volumetric capacities \(C^-\) and \(C^+\) and Backward Euler in time,
\begin{align}
  C^-\,\Delta x\,(\Phi_1^{n+1}-\Phi_1^n)
  &= \Delta t\,(q_{3/2}-q_{1/2}),\label{eq:1D-1-Phi}\\
  C^-\,\ell^-\,\bigl(\Phi_{2}^{-,n+1}-\Phi_{2}^{-,n}\bigr)
  &= \Delta t\,(q_{\Gamma}^- - q_{3/2}),\label{eq:1D-2-Phi}\\
  C^+\,\ell^+\,\bigl(\Phi_{2}^{+,n+1}-\Phi_{2}^{+,n}\bigr)
  &= \Delta t\,(q_{5/2}-q_{\Gamma}^+),\label{eq:1D-3-Phi}\\
  C^+\,\Delta x\,(\Phi_3^{n+1}-\Phi_3^n)
  &= \Delta t\,(q_{7/2}-q_{5/2}).\label{eq:1D-4-Phi}
\end{align}

\paragraph{Flux definitions}  
Using midpoint‐rule approximations,
\[
  q_{1/2} = a,
  \quad
  q_{3/2} = k^-\,\frac{\Phi_{2}^{-,n+1}-\Phi_1^{n+1}}{\tfrac{\Delta x}{2}+\ell^-},
  \quad
  q_{\Gamma}^- = k^-\,\frac{\Phi_{\Gamma}^{-,n+1} - \Phi_{2}^{-,n+1}}{\ell^-},
\]
\[
  q_{\Gamma}^+ = k^+\,\frac{\Phi_{2}^{+,n+1} - \Phi_{\Gamma}^{+,n+1}}{\ell^+},
  \quad
  q_{5/2} = k^+\,\frac{\Phi_3^{n+1}-\Phi_{2}^{+,n+1}}{\ell^+ + \tfrac{\Delta x}{2}},
  \quad
  q_{7/2} = b.
\]

\paragraph{Interface conditions}  
Continuity of the field and flux at \(x_\Gamma\) reads
\[
  \Phi_{\Gamma}^{-,n+1} - \Phi_{\Gamma}^{+,n+1} = 0,
  \qquad
  q_{\Gamma}^- - q_{\Gamma}^+ = 0,
\]
i.e.
\[
  k^-\,\frac{\Phi_{\Gamma}^- - \Phi_{2}^-}{\ell^-}
  = k^+\,\frac{\Phi_{2}^+ - \Phi_{\Gamma}^+}{\ell^+},
  \quad
  \Phi_{\Gamma}^- = \Phi_{\Gamma}^+.
\]

\paragraph{Uniform‐mesh limit}  
If \(\ell^-=\ell^+=\tfrac{\Delta x}{2}\) and \(k^-=k^+=k\), then
\[
  q_{3/2}=\frac{k}{\Delta x}(\Phi_{2}^- - \Phi_1),
  \quad
  q_{5/2}=\frac{k}{\Delta x}(\Phi_3 - \Phi_{2}^+),
\]
and \eqref{eq:1D-1-Phi}–\eqref{eq:1D-4-Phi} reduce to the standard second‐order finite‐volume scheme
\[
  C\,\Delta x\,(\Phi_i^{n+1}-\Phi_i^n)
  = \Delta t\,k\Bigl[
    \frac{\Phi_{i+1}^{n+1}-\Phi_i^{n+1}}{\Delta x}
    -\frac{\Phi_i^{n+1}-\Phi_{i-1}^{n+1}}{\Delta x}
  \Bigr],
\]
demonstrating consistency.

\paragraph{Assembled linear system}  
Collect the six unknowns at \(n+1\) into
\[
  \boldsymbol{\Phi}
  = [\,\Phi_1,\;\Phi_{2}^-,\;\Phi_{2}^+,\;\Phi_{\Gamma}^-,\;\Phi_{\Gamma}^+,\;\Phi_3\,]^{T}.
\]
Define
\[
  m_1 = C^-\,\Delta x,\quad
  m_2 = C^-\,\ell^-,\quad
  m_3 = C^+\,\ell^+,\quad
  m_4 = C^+\,\Delta x,
\]
\[
  \alpha_1 = \frac{k^-}{\tfrac{\Delta x}{2}+\ell^-},\quad
  \alpha_2 = \frac{k^-}{\ell^-},\quad
  \alpha_3 = \frac{k^+}{\ell^+},\quad
  \alpha_4 = \frac{k^+}{\ell^+ + \tfrac{\Delta x}{2}},
\]
and let \(a = q_{1/2}\), \(b = q_{7/2}\).  Then the discrete system
\(\mathbf{A}\,\boldsymbol{\Phi}=\mathbf{b}\) is
\[
\footnotesize
\mathbf{A} =
\begin{pmatrix}
m_1 + \Delta t\,\alpha_1 & -\Delta t\,\alpha_1 & 0 & 0 & 0 & 0 \\[6pt]
-\Delta t\,\alpha_1 & m_2 + \Delta t(\alpha_1+\alpha_2) & 0 & -\Delta t\,\alpha_2 & 0 & 0 \\[6pt]
0 & 0 & m_3 + \Delta t(\alpha_3+\alpha_4) & 0 & -\Delta t\,\alpha_3 & -\Delta t\,\alpha_4 \\[6pt]
0 & 0 & 0 & 1 & -1 & 0 \\[4pt]
0 & -\alpha_2 & -\alpha_3 & \alpha_2 & \alpha_3 & 0 \\[6pt]
0 & 0 & -\Delta t\,\alpha_4 & 0 & 0 & m_4 + \Delta t\,\alpha_4
\end{pmatrix}, 
\quad
\mathbf{b} =
\begin{pmatrix}
m_1\,\Phi_1^n + \Delta t\,a \\[4pt]
m_2\,\Phi_{2}^{-,n} \\[4pt]
m_3\,\Phi_{2}^{+,n} \\[4pt]
0 \\[4pt]
0 \\[4pt]
m_4\,\Phi_3^n + \Delta t\,b
\end{pmatrix}.
\]
Rows 1-3 and 6 enforce the finite‐volume balances \eqref{eq:1D-1-Phi}-\eqref{eq:1D-4-Phi}, row 4 enforces \(\Phi_{\Gamma}^- = \Phi_{\Gamma}^+\), and row 5 enforces \(q_{\Gamma}^- = q_{\Gamma}^+\).  This sparse system may be solved via a direct factorization or an iterative Krylov method with ILU preconditioning.

%% file: references.bib
@article{verzicco_immersed_2023,
	title = {Immersed {Boundary} {Methods}: {Historical} {Perspective} and {Future} {Outlook}},
	volume = {55},
	copyright = {http://creativecommons.org/licenses/by/4.0/},
	issn = {0066-4189, 1545-4479},
	shorttitle = {Immersed {Boundary} {Methods}},
	url = {https://www.annualreviews.org/doi/10.1146/annurev-fluid-120720-022129},
	doi = {10.1146/annurev-fluid-120720-022129},
	language = {en},
	number = {1},
	urldate = {2026-01-06},
	journal = {Annual Review of Fluid Mechanics},
	author = {Verzicco, Roberto},
	month = jan,
	year = {2023},
	pages = {129--155},
}

@article{peskin_flow_1972,
	title = {Flow patterns around heart valves: {A} numerical method},
	volume = {10},
	copyright = {https://www.elsevier.com/tdm/userlicense/1.0/},
	issn = {00219991},
	shorttitle = {Flow patterns around heart valves},
	url = {https://linkinghub.elsevier.com/retrieve/pii/0021999172900654},
	doi = {10.1016/0021-9991(72)90065-4},
	language = {en},
	number = {2},
	urldate = {2025-12-23},
	journal = {Journal of Computational Physics},
	author = {Peskin, Charles S},
	month = oct,
	year = {1972},
	pages = {252--271},
}

@article{johansen_cartesian_1998,
	title = {A {Cartesian} {Grid} {Embedded} {Boundary} {Method} for {Poisson}'s {Equation} on {Irregular} {Domains}},
	volume = {147},
	copyright = {https://www.elsevier.com/tdm/userlicense/1.0/},
	issn = {00219991},
	url = {https://linkinghub.elsevier.com/retrieve/pii/S0021999198959654},
	doi = {10.1006/jcph.1998.5965},
	language = {en},
	number = {1},
	urldate = {2025-06-26},
	journal = {Journal of Computational Physics},
	author = {Johansen, Hans and Colella, Phillip},
	month = nov,
	year = {1998},
	pages = {60--85},
}

@article{brown_diffusion_1965,
	title = {Diffusion of {Heat} {From} a {Sphere} to a {Surrounding} {Medium}},
	volume = {18},
	copyright = {https://doi.org/10.1071/journalslicense},
	issn = {0004-9506, 1446-5582},
	url = {https://connectsci.au/ph/article/18/5/483/197957/Diffusion-of-Heat-From-a-Sphere-to-a-Surrounding},
	doi = {10.1071/PH650483},
	language = {en},
	number = {5},
	urldate = {2025-12-17},
	journal = {Australian Journal of Physics},
	author = {Brown, A},
	month = oct,
	year = {1965},
	pages = {483--490},
}

@article{schwartz_cartesian_2006,
	title = {A {Cartesian} grid embedded boundary method for the heat equation and {Poisson}’s equation in three dimensions},
	volume = {211},
	copyright = {https://www.elsevier.com/tdm/userlicense/1.0/},
	issn = {00219991},
	url = {https://linkinghub.elsevier.com/retrieve/pii/S002199910500286X},
	doi = {10.1016/j.jcp.2005.06.010},
	language = {en},
	number = {2},
	urldate = {2025-11-13},
	journal = {Journal of Computational Physics},
	author = {Schwartz, Peter and Barad, Michael and Colella, Phillip and Ligocki, Terry},
	month = jan,
	year = {2006},
	pages = {531--550},
}

@article{verstappen_symmetry-preserving_2004,
	title = {A symmetry-preserving {Cartesian} grid method for computing a viscous flow past a circular cylinder},
	volume = {333},
	issn = {1631-0721, 1873-7234},
	url = {https://comptes-rendus.academie-sciences.fr/mecanique/articles/10.1016/j.crme.2004.09.021/},
	doi = {10.1016/j.crme.2004.09.021},
	language = {en},
	number = {1},
	urldate = {2025-12-17},
	journal = {Comptes Rendus. Mécanique},
	author = {Verstappen, Roel and Dröge, Marc},
	month = dec,
	year = {2004},
	pages = {51--57},
}

@book{Ishii2011,
	address = {New York, NY},
	title = {Thermo-{Fluid} {Dynamics} of {Two}-{Phase} {Flow}},
	copyright = {https://www.springernature.com/gp/researchers/text-and-data-mining},
	isbn = {978-1-4419-7984-1 978-1-4419-7985-8},
	url = {https://link.springer.com/10.1007/978-1-4419-7985-8},
	language = {en},
	urldate = {2025-11-03},
	publisher = {Springer},
	author = {Ishii, M. and Hibiki, T.},
	year = {2011},
	doi = {10.1007/978-1-4419-7985-8},
}

@article{pollack_kapitza_1969,
	title = {Kapitza {Resistance}},
	volume = {41},
	copyright = {http://link.aps.org/licenses/aps-default-license},
	issn = {0034-6861},
	url = {https://link.aps.org/doi/10.1103/RevModPhys.41.48},
	doi = {10.1103/RevModPhys.41.48},
	language = {en},
	number = {1},
	urldate = {2025-10-13},
	journal = {Reviews of Modern Physics},
	author = {Pollack, Gerald L.},
	month = jan,
	year = {1969},
	pages = {48--81},
}

@article{Sander2023,
	title = {Compilation of {Henry}'s law constants (version 5.0.0) for water as solvent},
	volume = {23},
	issn = {1680-7316},
	url = {https://acp.copernicus.org/articles/23/10901/2023/},
	doi = {10.5194/acp-23-10901-2023},
	language = {English},
	number = {19},
	urldate = {2025-08-28},
	journal = {Atmospheric Chemistry and Physics},
	author = {Sander, R.},
	month = oct,
	year = {2023},
	pages = {10901--12440},
}

@misc{santilli_geos_2023,
	title = {{GEOS} ({Geometry} {Engine}, {Open} {Source})},
	copyright = {GNU Lesser General Public License v2.1 only},
	url = {https://zenodo.org/doi/10.5281/zenodo.11396894},
	urldate = {2025-07-04},
	publisher = {Zenodo},
	author = {Santilli, Sandro and Ramsey, Paul and Baston, Daniel and Łoskot, Mateusz and Davis, Martin and Obe, Regina and Savage, Charlie and Bychkov, Yury and Rouault, Even and Toews, Michael},
	month = nov,
	year = {2023},
	doi = {10.5281/ZENODO.11396894},
}

@article{saye_high-order_2015,
	title = {High-{Order} {Quadrature} {Methods} for {Implicitly} {Defined} {Surfaces} and {Volumes} in {Hyperrectangles}},
	volume = {37},
	issn = {1064-8275, 1095-7197},
	url = {http://epubs.siam.org/doi/10.1137/140966290},
	doi = {10.1137/140966290},
	language = {en},
	number = {2},
	urldate = {2025-07-04},
	journal = {SIAM Journal on Scientific Computing},
	author = {Saye, R. I.},
	month = jan,
	year = {2015},
	pages = {A993--A1019},
}

@article{remmerswaal_sharp_2022,
	title = {A sharp, structure preserving two-velocity model for two-phase flow},
	url = {https://rgdoi.net/10.13140/RG.2.2.14099.55841},
	doi = {10.13140/RG.2.2.14099.55841},
	language = {en},
	urldate = {2025-07-04},
	author = {Remmerswaal, Ronald A and Veldman, Arthur E P},
	year = {2022},
}

@inproceedings{laidlaw_constructive_1986,
	title = {Constructive solid geometry for polyhedral objects},
	isbn = {9780897911962},
	url = {https://dl.acm.org/doi/10.1145/15922.15904},
	doi = {10.1145/15922.15904},
	language = {en},
	urldate = {2025-07-04},
	booktitle = {Proceedings of the 13th annual conference on {Computer} graphics and interactive techniques},
	publisher = {ACM},
	author = {Laidlaw, David H. and Trumbore, W. Benjamin and Hughes, John F.},
	month = aug,
	year = {1986},
	pages = {161--170},
}

@article{gabbard_high-order_2024,
	title = {A high-order finite difference method for moving immersed domain boundaries and material interfaces},
	volume = {507},
	issn = {00219991},
	url = {https://linkinghub.elsevier.com/retrieve/pii/S0021999124002286},
	doi = {10.1016/j.jcp.2024.112979},
	language = {en},
	urldate = {2025-07-04},
	journal = {Journal of Computational Physics},
	author = {Gabbard, James and Van Rees, Wim M.},
	month = jun,
	year = {2024},
	pages = {112979},
}

@article{gibou_second-order-accurate_2002,
	title = {A {Second}-{Order}-{Accurate} {Symmetric} {Discretization} of the {Poisson} {Equation} on {Irregular} {Domains}},
	volume = {176},
	copyright = {https://www.elsevier.com/tdm/userlicense/1.0/},
	issn = {00219991},
	url = {https://linkinghub.elsevier.com/retrieve/pii/S0021999101969773},
	doi = {10.1006/jcph.2001.6977},
	language = {en},
	number = {1},
	urldate = {2025-07-04},
	journal = {Journal of Computational Physics},
	author = {Gibou, Frederic and Fedkiw, Ronald P. and Cheng, Li-Tien and Kang, Myungjoo},
	month = feb,
	year = {2002},
	pages = {205--227},
}

@article{cheny_ls-stag_2010,
	title = {The {LS}-{STAG} method: {A} new immersed boundary/level-set method for the computation of incompressible viscous flows in complex moving geometries with good conservation properties},
	volume = {229},
	copyright = {https://www.elsevier.com/tdm/userlicense/1.0/},
	issn = {00219991},
	shorttitle = {The {LS}-{STAG} method},
	url = {https://linkinghub.elsevier.com/retrieve/pii/S002199910900552X},
	doi = {10.1016/j.jcp.2009.10.007},
	language = {en},
	number = {4},
	urldate = {2025-07-03},
	journal = {Journal of Computational Physics},
	author = {Cheny, Yoann and Botella, Olivier},
	month = feb,
	year = {2010},
	pages = {1043--1076},
}

@book{carslaw_conduction_1980,
	address = {Oxford},
	edition = {2. ed., repr},
	title = {Conduction of heat in solids},
	isbn = {9780198533030},
	language = {eng},
	publisher = {Clarendon Pr},
	author = {Carslaw, Horatio S. and Jaeger, J. C.},
	year = {1980},
}

@article{popinet_accurate_2009,
	title = {An accurate adaptive solver for surface-tension-driven interfacial flows},
	volume = {228},
	copyright = {https://www.elsevier.com/tdm/userlicense/1.0/},
	issn = {00219991},
	url = {https://linkinghub.elsevier.com/retrieve/pii/S002199910900240X},
	doi = {10.1016/j.jcp.2009.04.042},
	language = {en},
	number = {16},
	urldate = {2025-06-27},
	journal = {Journal of Computational Physics},
	author = {Popinet, Stéphane},
	month = sep,
	year = {2009},
	pages = {5838--5866},
}

@article{hirt_arbitrary_1974,
	title = {An arbitrary {Lagrangian}-{Eulerian} computing method for all flow speeds},
	volume = {14},
	copyright = {https://www.elsevier.com/tdm/userlicense/1.0/},
	issn = {00219991},
	url = {https://linkinghub.elsevier.com/retrieve/pii/0021999174900515},
	doi = {10.1016/0021-9991(74)90051-5},
	language = {en},
	number = {3},
	urldate = {2025-06-26},
	journal = {Journal of Computational Physics},
	author = {Hirt, C.W and Amsden, A.A and Cook, J.L},
	month = mar,
	year = {1974},
	pages = {227--253},
}

@article{mirtich_fast_1996,
	title = {Fast and {Accurate} {Computation} of {Polyhedral} {Mass} {Properties}},
	volume = {1},
	issn = {1086-7651},
	url = {http://www.tandfonline.com/doi/abs/10.1080/10867651.1996.10487458},
	doi = {10.1080/10867651.1996.10487458},
	language = {en},
	number = {2},
	urldate = {2025-06-26},
	journal = {Journal of Graphics Tools},
	author = {Mirtich, Brian},
	month = jan,
	year = {1996},
	pages = {31--50},
}

@article{hirt_volume_1981,
	title = {Volume of fluid ({VOF}) method for the dynamics of free boundaries},
	volume = {39},
	copyright = {https://www.elsevier.com/tdm/userlicense/1.0/},
	issn = {00219991},
	url = {https://linkinghub.elsevier.com/retrieve/pii/0021999181901455},
	doi = {10.1016/0021-9991(81)90145-5},
	language = {en},
	number = {1},
	urldate = {2025-06-26},
	journal = {Journal of Computational Physics},
	author = {Hirt, C.W and Nichols, B.D},
	month = jan,
	year = {1981},
	pages = {201--225},
}

@article{osher_fronts_1988,
	title = {Fronts propagating with curvature-dependent speed: {Algorithms} based on {Hamilton}-{Jacobi} formulations},
	volume = {79},
	copyright = {https://www.elsevier.com/tdm/userlicense/1.0/},
	issn = {00219991},
	shorttitle = {Fronts propagating with curvature-dependent speed},
	url = {https://linkinghub.elsevier.com/retrieve/pii/0021999188900022},
	doi = {10.1016/0021-9991(88)90002-2},
	language = {en},
	number = {1},
	urldate = {2025-06-26},
	journal = {Journal of Computational Physics},
	author = {Osher, Stanley and Sethian, James A},
	month = nov,
	year = {1988},
	pages = {12--49},
}

@book{cussler_diffusion_2009,
	edition = {3},
	title = {Diffusion: {Mass} {Transfer} in {Fluid} {Systems}},
	copyright = {https://www.cambridge.org/core/terms},
	isbn = {9780521871211 9780511805134},
	shorttitle = {Diffusion},
	url = {https://www.cambridge.org/core/product/identifier/9780511805134/type/book},
	urldate = {2025-06-26},
	publisher = {Cambridge University Press},
	author = {Cussler, E. L.},
	month = jan,
	year = {2009},
	doi = {10.1017/CBO9780511805134},
}

@article{mccorquodale_high-order_2011,
	title = {A high-order finite-volume method for conservation laws on locally refined grids},
	volume = {6},
	issn = {2157-5452, 1559-3940},
	url = {http://msp.org/camcos/2011/6-1/p01.xhtml},
	doi = {10.2140/camcos.2011.6.1},
	language = {en},
	number = {1},
	urldate = {2025-06-26},
	journal = {Communications in Applied Mathematics and Computational Science},
	author = {McCorquodale, Peter and Colella, Phillip},
	month = mar,
	year = {2011},
	pages = {1--25},
}

@article{calhoun_cartesian_2000,
	title = {A {Cartesian} {Grid} {Finite}-{Volume} {Method} for the {Advection}-{Diffusion} {Equation} in {Irregular} {Geometries}},
	volume = {157},
	copyright = {https://www.elsevier.com/tdm/userlicense/1.0/},
	issn = {00219991},
	url = {https://linkinghub.elsevier.com/retrieve/pii/S0021999199963696},
	doi = {10.1006/jcph.1999.6369},
	language = {en},
	number = {1},
	urldate = {2025-06-26},
	journal = {Journal of Computational Physics},
	author = {Calhoun, Donna and LeVeque, Randall J.},
	month = jan,
	year = {2000},
	pages = {143--180},
}

@book{bird_transport_2007,
	address = {New York},
	edition = {Revised ed},
	title = {Transport phenomena},
	isbn = {9780470115398},
	language = {eng},
	publisher = {Wiley},
	author = {Bird, Robert Byron and Stewart, Warren E. and Lightfoot, Edwin N.},
	year = {2007},
}

@book{polyanin_handbook_2001,
	edition = {0},
	title = {Handbook of {Linear} {Partial} {Differential} {Equations} for {Engineers} and {Scientists}},
	isbn = {9780429140938},
	url = {https://www.taylorfrancis.com/books/9781420035322},
	language = {en},
	urldate = {2025-06-24},
	publisher = {Chapman and Hall/CRC},
	author = {Polyanin, Andrei D.},
	month = nov,
	year = {2001},
	doi = {10.1201/9781420035322},
}

@article{limare_hybrid_2023,
	title = {A hybrid level-set / embedded boundary method applied to solidification-melt problems},
	volume = {474},
	issn = {00219991},
	url = {https://linkinghub.elsevier.com/retrieve/pii/S0021999122008920},
	doi = {10.1016/j.jcp.2022.111829},
	language = {en},
	urldate = {2025-05-23},
	journal = {Journal of Computational Physics},
	author = {Limare, A. and Popinet, S. and Josserand, C. and Xue, Z. and Ghigo, A.},
	month = feb,
	year = {2023},
	pages = {111829},
}

@article{li_numerical_2019,
	title = {Numerical {Simulation} of {Stefan} {Problem} {Coupled} with {Mass} {Transport} in a {Binary} {System} {Through} {XFEM}/{Level} {Set} {Method}},
	volume = {78},
	issn = {0885-7474, 1573-7691},
	url = {http://link.springer.com/10.1007/s10915-018-0759-x},
	doi = {10.1007/s10915-018-0759-x},
	language = {en},
	number = {1},
	urldate = {2025-05-22},
	journal = {Journal of Scientific Computing},
	author = {Li, Min and Chaouki, Hicham and Robert, Jean-Loup and Ziegler, Donald and Fafard, Mario},
	month = jan,
	year = {2019},
	pages = {145--166},
}

@article{chierici_optimized_2022,
	title = {An optimized {Vofi} library to initialize the volume fraction field},
	volume = {281},
	issn = {00104655},
	url = {https://linkinghub.elsevier.com/retrieve/pii/S0010465522002259},
	doi = {10.1016/j.cpc.2022.108506},
	language = {en},
	urldate = {2025-05-20},
	journal = {Computer Physics Communications},
	author = {Chierici, A. and Chirco, L. and Le Chenadec, V. and Scardovelli, R. and Yecko, Ph. and Zaleski, S.},
	month = dec,
	year = {2022},
	pages = {108506},
}

@misc{rodriguez_conservative_2022,
	title = {A {Conservative} {Cartesian} {Cut} {Cell} {Method} for the {Solution} of the {Incompressible} {Navier}-{Stokes} {Equations} on {Staggered} {Meshes}},
	copyright = {arXiv.org perpetual, non-exclusive license},
	url = {https://arxiv.org/abs/2211.10698},
	doi = {10.48550/ARXIV.2211.10698},
	urldate = {2025-05-19},
	publisher = {arXiv},
	author = {Rodríguez, Alejandro Quirós and Fullana, Tomas and Chenadec, Vincent Le and Sayadi, Taraneh},
	year = {2022},
	note = {Version Number: 1},
}
